\documentclass[12pt, a4paper]{article}
\usepackage{amssymb,amsmath,amsfonts,amsthm}
\usepackage{setspace}
\usepackage[DIV14]{typearea}
\usepackage[stretch=10]{microtype}
\usepackage{hyperref}

\DeclareMathOperator{\BV}{BV}
\numberwithin{equation}{section}

\newcommand{\E}{\mathbb{E}}

\usepackage{dsfont}

\theoremstyle{plain}
\newtheorem{thm}{Theorem}[section]
\newtheorem{lma}[thm]{Lemma}
\newtheorem{rmk}[thm]{Remark}
\newtheorem{prop}[thm]{Proposition}
\newtheorem{cor}[thm]{Corollary}

\newtheorem{defn}[thm]{Definition}

\newtheorem{claim}[thm]{Claim}

\theoremstyle{definition}

\begin{document}

\title{A quenched weak invariance principle }
\author{{\small J\'{e}r\^{o}me Dedecker$^a$, Florence Merlev\`ede$^b$ and Magda
Peligrad$^c$ \thanks{%
Supported in part by a Charles Phelps Taft Memorial Fund grant,
the NSA grant H98230-11-1-0135, and the NSF grant DMS-1208237.}}}

\date{}
\maketitle

\noindent $^a$ Universit\'e Paris Descartes, Sorbonne Paris Cit\'e, Laboratoire MAP5
and CNRS UMR 8145.  \\
$^b$ Universit\'e Paris Est, LAMA (UMR 8050), UPEMLV, CNRS, UPEC. \\
$^c$ University of Cincinnati, Department of Mathematical Sciences.

\begin{center}
\textbf{Abstract}
\end{center}

\noindent In this paper we study the almost sure conditional central limit theorem in
its functional form for a  class of random variables satisfying a
projective criterion. Applications to strongly mixing processes and non irreducible
Markov chains are given. The proofs are based on the normal
approximation of double indexed martingale-like sequences, an approach which has
interest in itself.

\vskip5pt

\begin{center}
\textbf{R\'esum\'e}
\end{center}

\noindent Dans cet article, nous \'etudions le th\'eor\`eme central limite conditionnel presque s\^ur, ainsi que sa forme fonctionnelle, pour des suites stationnaires de variables al\'eatoires r\'eelles satisfaisant une condition de type projectif. Nous donnons des applications de ces r\'esultats aux processus fortement m\'elangeants ainsi qu'\`a des cha\^ines de Markov non irr\'eductibles. Les preuves sont essentiellement bas\'ees sur une approximation normale de suites doublement index\'ees de variables al\'eatoires de type martingale.

\vskip5pt


\vskip5pt


\section{Introduction}

Let $(\xi_{i})_{i \geq 0}$ be a Markov chain admitting an invariant probability
$\pi$. Let
$f$ be a real-valued function such that
$\pi(f^2)< \infty$ and $\pi(f)=0$, and let $S_n=f(\xi_1)+ \cdots + f(\xi_n)$.
If the central limit theorem (CLT)  holds for $n^{-1/2}S_n$ starting form the initial
distribution $\pi$, an interesting question is to
know whether it remains true for another initial distribution $\nu$.
 Maxwell and Woodroofe (2000) have given a projective criterion
 under which $S_n$ satisfies the so-called conditional CLT, which implies that the CLT holds
for any initial distribution having a bounded density with respect to $\pi$. Necessary and sufficient conditions
for the conditional CLT are given in  Dedecker and Merlev\`ede (2002), and
Wu and Woodroofe (2004).

The question is more delicate
if $\nu$ is a Dirac mass at point $x$. One says that the CLT is quenched
if it holds for almost every starting point with respect to $\pi$. The quenched CLT implies the central limit theorem for the chain starting from an invariant probability measure $\pi$, referred as annealed CLT. The same terminologies
are used for the functional central limit theorem (FCLT). 
  For aperiodic Harris recurrent Markov chains, the quenched CLT question is solved by using 
Proposition 18.1.2 in Meyn and Tweedie (1993). More precisely, for an aperiodic Harris recurrent Markov chain, if the CLT holds for the initial
distribution $\pi$, then it holds for any initial distribution, and hence for any
starting point $x$ (see Proposition 3.1 in Chen (1999) and its proof). In the non irreducible setting, the situation is not so clear. For instance, an example of a Markov chain with 
normal transition operator satisfying the annealed CLT but not the quenched is given at the end of Section 3 in Derriennic and Lin (2001). 

This question of the quenched CLT can be formulated in the more general context
of stationary sequences: it means that, on a set of measure one,  the central limit theorem holds
when replacing the usual expectation by the conditional expectation with respect to
the past $\sigma$-algebra. Some examples of stationary processes satisfying the CLT
 but not the quenched CLT can be found in  Voln\'{y} and Woodroofe (2010a).

The first general results on the quenched CLT and FCLT are given in Borodin and
Ibragimov (1994): in the Markov chain setting, it says that the FCLT holds if there is a solution
in ${\mathbb L}^2(\pi)$ to the Poisson equation (see Gordin and Lifschitz (1978)); in
a general setting it means that the FCLT is true under  Gordin's condition (1969).
This result has been improved by
Derriennic and Lin (2001, 2003),
Zhao and Woodroofe (2008), Cuny (2011), Cuny and
Peligrad (2012), Cuny and Voln\'y (2012), Voln\'{y} and Woodroofe (2010b) and Merlev\`{e}de \textit{et
al.} (2012). In a recent paper, Cuny and Merlev\`{e}de (2012) have proved that
the FCLT is quenched under the condition of Maxwell and Woodroofe (2000).

All the papers cited  above use a martingale approximation in
${\mathbb L}^2$. Consequently, the  projective condition obtained up to now are
always  expressed in terms
of ${\mathbb L}^2$ norms of conditional expectations.
In this paper, we prove the quenched FCLT under a projective condition
involving ${\mathbb L}^1$-norms, in the spirit of Gordin (1973).
As a consequence, we obtain that the FCLT of Doukhan \textit{et
al.} (1994)
for strongly mixing sequences is quenched. Note that Doukhan \textit{et al} (1994)
have shown that their condition is optimal in some sense for the usual FCLT,
so it is also sharp for the quenched FCLT. In Section 3.1, we study the example
of the non irreducible Markov chain associated to an intermittent map. Once again, we shall see through this example that our condition is essentially optimal.

Our main result, Theorem \ref{mainth} below,  is a consequence of the more general Proposition \ref{mainprop}, where the conditions are expressed in terms
of conditional expectations of partial sums.
The proof of this proposition is done via a blocking argument followed by a two
step martingale decomposition. We start with a finite number of consecutive
blocks of random variables. The sum in blocks are approximated by
martingales. This decomposition introduces the need of studying the normal
approximation for a family of double indexed martingales. This approximation
has interest in itself and is presented in Section \ref{NA}.

\section{Results}\label{results}

Let $(\Omega,\mathcal{A}, {\mathbb{P}})$ be a probability space, and $%
T:\Omega \mapsto \Omega$ be a bijective bimeasurable transformation
preserving the probability ${\mathbb{P}}$. An element $A$  is said to be
invariant if $T(A)=A$. We denote by $\mathcal{I}$ the $\sigma$-algebra of
all invariant sets. The probability ${\mathbb{P}}$ is ergodic if each
element of $\mathcal{I}$ has measure 0 or 1.

Let $\mathcal{F}_0$ be a $\sigma$-algebra of $\mathcal{A}$ satisfying $%
\mathcal{F}_0 \subseteq T^{-1}(\mathcal{F}_0)$ and define the nondecreasing
filtration $(\mathcal{F}_i)_{i \in {\mathbb{Z}}}$ by $\mathcal{F}_i= T^{-i}(%
\mathcal{F}_0)$. We assume that there exists a regular version $P_{T|{%
\mathcal{F}}_0}$ of $T$ given ${\mathcal{F}}_0$, and for any integrable
random variable $f$ from $\Omega$ to ${\mathbb{R}}$ we write $ K(f)=P_{T|{%
\mathcal{F}}_0}(f) $. Since ${\mathbb{P}}$ is invariant by $T$, for any
integer $k$, a regular version $P_{T|{\mathcal{F}}_{k}}$ of $T$ given ${%
\mathcal{F}}_{k}$ is then obtained $via$ $P_{T|{\mathcal{F}}%
_{k}}(f)=K(f\circ T^{-k}) \circ T^{k}$. In the sequel, all the conditional
expectations with respect to ${\mathcal{F}}_k$ are
obtained through these conditional probabilities. More precisely, we shall use the following notations: 
$$
\mathbb{E}_0(X) :=\mathbb{E} (X |{\mathcal{F}}_0)= K(X \circ T^{-1}) \, \text{ and } \, \mathbb{E}_k(X) := \mathbb{E} (X |{\mathcal{F}}_k) = K(X \circ T^{-k-1}) \circ T^k \, .
$$
With these notations, ${\mathbb{E}}(f \circ T^2 |{\mathcal{F}}_0)=\mathbb{E}(K(f)\circ T|{\mathcal{F}}_0)=K^2(f)$, and more
generally, for any positive integer $\ell$, ${\mathbb{E}}(f \circ T^\ell |{%
\mathcal{F}}_0)=K^\ell(f)$. 

Let $X_0$ be an $\mathcal{F}_0$-measurable, square integrable and centered
random variable. Define the sequence $\mathbf{X}=(X_i)_{i \in \mathbb{Z}}$
by $X_i=X_0\circ T^i$. Let $S_n=X_1+ \cdots + X_n$, and define the Donsker
process $W_n$ by $W_n(t)=n^{-1/2}(S_{[nt]}+(nt-[nt])X_{[nt]+1})$.

Let $\mathcal{H}^{*}$ the space of continuous functions $\varphi$ from $%
(C([0,1]), \|.\|_{\infty})$ to ${\mathbb{R}}$ such that $x\rightarrow
|(1+\|x\|_{\infty}^2)^{-1}\varphi(x)|$ is bounded. Our main result is the
following:

\begin{thm}
\label{mainth} Assume that
\begin{equation}
\sum_{k\geq 0}\Vert X_{0}{\mathbb{E}}_{0}(X_{k})\Vert _{1}<\infty \,.
\label{cond}
\end{equation}%
then the series
\begin{equation}
\eta ={\mathbb{E}}(X_{0}^{2}|\mathcal{I})+2\sum_{k>0}{\mathbb{E}}(X_{0}X_{k}|%
\mathcal{I})  \label{eta}
\end{equation}%
converges almost surely and in ${\mathbb{L}}^{1}$. Moreover, on a set of probability one, for any $\varphi $ in $\mathcal{H}%
^{\ast }$, 
\begin{equation}
\lim_{n\rightarrow \infty }{\mathbb{E}}_{0}(\varphi (W_{n}))=\int \varphi (z
\sqrt{\eta })W(dz) \, ,
\label{L1ps}
\end{equation}%
where $W$ is the distribution of a standard Wiener process. The convergence in \eqref{L1ps} also holds in ${\mathbb{L}}^{1}$. 
\end{thm}

Note that the ${\mathbb{L}}^{1}$-convergence in (\ref{L1ps}) has been proved
in Dedecker and Merlev\`{e}de (2002). In this paper, we shall prove the
almost sure convergence. Various classes of examples satisfying (\ref{cond})
can be found in Dedecker and Rio (2000).

This result has an interesting interpretation in the terminology of additive
functionals of Markov chains.
Let $(\xi _{n})_{n\geq 0}$ be a Markov chain with values in a Polish space $S$,
so that there exists a regular transition probability $P_{\xi_1|\xi_0=x}$. Let  $P$
be the transition kernel defined by $%
P(f)(x )=P_{\xi_1|\xi_0=x}(f)$ for any  bounded measurable function $f$
from $S$ to ${\mathbb R}$, and assume that
there exists an invariant probability $\pi$ for this transition kernel,
that is a probability measure on $S$ such that $\pi(f)=\pi(P(f))$ for any
bounded measurable function $f$ from $S$ to ${\mathbb R}$.
Let then $\mathbb{L}_{0}^{2}(\pi )$ be
the set of functions from $S$ to ${\mathbb R}$ such that $\pi(f^2) <\infty $ and $\pi(f)=0$.
For $f\in \mathbb{L}_{0}^{2}(\pi )${\ define $%
X_{i}=f(\xi _{i})$.
Notice that any
stationary sequence $(Y_{k})_{k\in \mathbb{Z}}$ can be viewed as a function
of a Markov process $\xi _{k}=(Y_{i};i\leq k),$ for the function $g(\xi
_{k})=Y_{k}$.

In this setting the condition (\ref{cond}) is $\sum_{k\geq 0}
\pi(|fP^{k}(f)|) <\infty .$
Also, the random variable $\eta$ defined in Theorem \ref{mainth}
 is the limit
almost surely and in ${\mathbb L}^1$ of $n^{-1}{\mathbb E}(S_n^2|\xi_0)$,
in such a way that $\eta=\bar\eta(\xi_0)$. By stationarity, it is also the limit in ${\mathbb L}^1$ of the sequence
$n^{-1}{\mathbb E}((X_2+ \cdots + X_{n+1})^2|\xi_1)$, so that
$\bar \eta(\xi_0)=\bar \eta(\xi_1)$ almost surely. Consequently $\bar \eta$ is an 
harmonic function for $P$ in the sense that  $\pi$-almost surely $P(\bar \eta)=\bar \eta$. 

In the context of Markov chain the conclusion of Theorem \ref{mainth} is also known under
the terminology of FCLT started at a point. To rephrase it, let $\mathbb{P}%
^{x}$ be the probability associated to the Markov chain started from $x$ and
let $\mathbb{E}^{x}$ be the corresponding expectation. Then, for $\pi$-almost every $x \in S$, for any $%
\varphi $ in ${\mathcal H}^*$,
\begin{equation*}
\lim_{n \rightarrow \infty} \mathbb{E}^{x}(\varphi (W_{n}))= \int \varphi (z\sqrt{\bar \eta(x)
})W(dz) \, . 
\end{equation*}
Moreover, 
\begin{equation*}
\lim_{n \rightarrow \infty} \int \Big|\mathbb{E}^{x}(\varphi (W_{n})) - \int \varphi (z\sqrt{\bar \eta(x)
})W(dz)\Big|\pi(dx)=0\, .
\end{equation*}

We mention
that in  Theorem \ref{mainth} no assumption of irreducibility nor of
aperiodicity is imposed. Under the additional assumptions that the Markov
chain is irreducible, aperiodic and positively recurrent, Chen (1999) showed
that the CLT holds for the stationary Markov chain under the condition
$\sum_{k\geq 0}\pi( fP^{k}(f))$ is convergent, and the quenched CLT holds
under the same condition by applying his Proposition 3.1.

\begin{rmk} \label{GordinL1} Let us present an alternative condition to the criterion \eqref{cond} in case where $T$ is ergodic. We do not require here $X_0$ to be in ${\mathbb L}^2$ but only in ${\mathbb L}^1$. The so-called Gordin criterion in ${\mathbb L}^1$ is:
\begin{equation} \label{condGordinL1}
\sup_{n \in {\mathbb N}} \Vert {\mathbb E}_0 (S_n) \Vert_1 < \infty \ \text{ and } \ \liminf_{n \rightarrow \infty} \frac{{\mathbb E} |S_n|}{ \sqrt n} < \infty \, .
\end{equation}
By Esseen and Janson (1985), it is known that  \eqref{condGordinL1} is equivalent to the following ${\mathbb L}^1$-coboundary decomposition:
\begin{equation} \label{coboundaryL1}
X_0 = m_0 +z_0 -z_0 \circ T \, ,
\end{equation}
where $z_0 \in {\mathbb L}^1$ and $m_0$ is a ${\mathcal F}_0$-measurable random variable in $ {\mathbb L}^2$ such that $\E_{-1} (m_0) = 0$ almost surely. Therefore, the criterion \eqref{GordinL1} leads to the annealed CLT. Note that one can easily prove that the condition \eqref{condalpha} of the next section also implies \eqref{condGordinL1}. However, the condition \eqref{condGordinL1} is not sufficient to get  the annealed FCLT (see Voln\'{y} and Samek (2000)). In addition, from Corollary 2 in Voln\'{y} and Woodroofe (2010,b), it follows that \eqref{condGordinL1} is not sufficient to get the quenched CLT either. In Proposition \ref{propDV} of Section \ref{sectionremark}, we shall provide  an example of stationary process for which \eqref{cond} holds but \eqref{condGordinL1} fails.

\end{rmk}

\section{Applications} \label{sectionappli}

\quad As a consequence of Theorem \ref{mainth}, we obtain the following
corollary for a class of weakly dependent sequences. We first need some definitions.

\begin{defn}
\label{defalphaweak} For a sequence $\mathbf{Y}=(Y_i)_{i \in {\mathbb{Z}}}$,
where $Y_i=Y_0 \circ T^i$ and $Y_0$ is an $\mathcal{F}_0$-measurable and
real-valued random variable, let for any $k \in {\mathbb{N}}$,
\begin{equation*}
\alpha_{{\mathbf{Y}}}(k) = \sup_{t \in {\mathbb{R}}} \big \| \mathbb{E} (
\mathbf{1}_{Y_k \leq t} | {\mathcal{F}}_0) - \mathbb{E} ( \mathbf{1}_{Y_k
\leq t}) \big \|_1 \, .
\end{equation*}
\end{defn}

\begin{defn}
\label{defalphastrong}
Recall that the strong mixing coefficient of Rosenblatt (1956) between two $\sigma$-algebras ${\mathcal F}$
and ${\mathcal G}$ is defined by $  \alpha({\mathcal F}, {\mathcal G})= \sup_{A \in {\mathcal F}, B \in {\mathcal G}}|{\mathbb P}(A \cap B)-{\mathbb P}(A){\mathbb P}(B)| $. 
For a strictly stationary sequence $(Y_i)_{i \in {\mathbb Z}}$ of real valued random variables, and the $\sigma$-algebra ${\mathcal F}_0=\sigma (Y_i, i \leq 0)$, define then
\begin{equation}\label{defalpharosen}
 \alpha(0) = 1 \text{ and } \alpha(k)= 2\alpha({\mathcal F}_0,\sigma(Y_k)) \text{ for $k>0$} \, .
\end{equation}
\end{defn}
Between the two above coefficients, the following relation holds: for any positive $k$, $\alpha_{ {\bf Y}}(k) \leq \alpha (k)$. In addition, the $\alpha$-dependent coefficient as defined in Definition \ref{defalphaweak} may be computed for instance for many Markov chains associated to dynamical systems that fail to be strongly mixing in the sense of Rosenblatt (see Section \ref{sectweakMC}).

\begin{defn}
\label{defquant}
A quantile function $Q$ is a function from $]0,1]$ to ${\mathbb R}_+$,
which is left-continuous and  non increasing. For any nonnegative random variable $Z$, we define the
 quantile function $Q_Z $ of $Z$  by $Q_Z (u) = \inf  \{ t \geq 0
: {\mathbb{P}} (|Z| >t ) \leq u \} $.
\end{defn}

\begin{defn}
\label{defclosedenv} Let $\mu$ be the probability distribution of a random
variable $X$. If $Q$ is an integrable quantile function, let $\mathrm{Mon}%
(Q, \mu)$ be the set of functions $g$ which are monotonic on some open
interval of ${\mathbb{R}}$ and null elsewhere and such that $Q_{|g(X)|} \leq
Q$. Let $\mathcal{F}(Q, \mu)$ be the closure in ${\mathbb{L}}^1(\mu)$ of the
set of functions which can be written as $\sum_{\ell=1}^{L} a_\ell f_\ell$,
where $\sum_{\ell=1}^{L} |a_\ell| \leq 1$ and $f_\ell$ belongs to $\mathrm{%
Mon}(Q, \mu)$.
\end{defn}

\begin{cor}
\label{coralpha} Let $Y_0$ be a real-valued random variable with law $P_{Y_0}$, and $Y_i=Y_0
\circ T^i$.  Let $Q$ be a quantile function such that
\begin{equation}
\sum_{k\geq 0}\int_{0}^{\alpha_{\mathbf{Y}} (k)}Q^{2}(u)du<\infty \, .
\label{condalpha}
\end{equation}%
Let $X_i = f(Y_i) - \mathbb{E} ( f(Y_i))$, where
 $f$ belongs to $\mathcal{F}(Q, P_{Y_0})$.
Then (\ref{cond}) is satisfied and consequently, the conclusion of Theorem \ref{mainth} holds.
\end{cor}

To prove that (\ref{condalpha}) implies (\ref{cond}), it suffices to apply
Proposition 5.3 with $m=q=1$ of Merlev\`ede and Rio (2012).

Notice that if $(\alpha(k))_{k \geq 0}$ is the usual sequence of strong
mixing coefficients of the stationary sequence $(X_i)_{i \in {\mathbb{Z}}}$
as defined in \eqref{defalpharosen}, then it follows from Corollary \ref%
{coralpha} that if
\begin{equation}
\sum_{k\geq 0}\int_{0}^{\alpha (k)}Q^{2}_{|X_0|}(u)du<\infty \, ,
\label{condalphafort}
\end{equation}
then the conclusion of Theorem \ref{mainth} holds.   Hence
the weak invariance principle of Doukhan \textit{et al.} (1994) is also
quenched. We refer to Theorem 2 in Doukhan \textit{et al.} (1994) and to
Bradley (1997) for a discussion on the optimality of the condition (\ref%
{condalphafort}).

\subsection{Application to functions of Markov chains associated to
intermittent maps}

\label{sectweakMC}

For $\gamma$ in $]0, 1[$, we consider the intermittent map $T_\gamma$ from $%
[0, 1]$ to $[0, 1]$, which is a modification of the Pomeau-Manneville map
(1980):
\begin{equation*}
T_\gamma(x)=
\begin{cases}
x(1+ 2^\gamma x^\gamma) \quad \text{ if $x \in [0, 1/2[$} \\
2x-1 \quad \quad \quad \ \ \text{if $x \in [1/2, 1]$} \, .%
\end{cases}
\end{equation*}
Recall that $T_{\gamma}$ is ergodic (and even mixing in the ergodic
theoretic sense) and that there exists a unique $T_\gamma$-invariant
probability measure  $\nu_\gamma$ on $[0, 1]$, which is absolutely
continuous with respect to the Lebesgue measure. We denote by $L_\gamma$ the
Perron-Frobenius operator of $T_\gamma$ with respect to $\nu_\gamma$. Recall
that for any bounded measurable functions $f$ and $g$, $
\nu_\gamma(f \cdot g\circ T_\gamma)=\nu_\gamma(L_\gamma(f) g) $. Let $(Y_i)_{i \geq 0}$ be a  Markov chain with transition Kernel $L_\gamma$  and invariant measure $%
\nu_\gamma$.

\begin{defn}
A function $H$ from ${\mathbb{R}}_+$ to $[0, 1]$ is a tail function if it is
non-increasing, right continuous, converges to zero at infinity, and $%
x\rightarrow x H(x)$ is integrable. If $\mu$ is a probability measure on $%
\mathbb{R}$ and $H$ is a tail function, let $\mathrm{Mon}^*(H, \mu)$ denote
the set of functions $f:\mathbb{R}\to \mathbb{R}$ which are monotonic on
some open interval and null elsewhere and such that $\mu(|f|>t)\leq H(t)$.
Let $\mathcal{F}^*(H, \mu)$ be the closure in $\mathbf{L}^1(\mu)$ of the set
of functions which can be written as $\sum_{\ell=1}^L a_\ell f_\ell$, where $%
\sum_{\ell=1}^L |a_\ell| \leq 1$ and $f_\ell\in \mathrm{Mon}^*(H, \mu)$.
\end{defn}

\begin{cor}
\label{ASmapB} Let $\gamma \in (0,1/2)$ and $(Y_i)_{i \geq 1}$ be a
stationary Markov chain with transition kernel $L_{\gamma}$ and invariant
measure $\nu_\gamma$. Let $H$ be a tail function such that
\begin{equation}  \label{lilcond}
\int_0^{\infty} x (H(x))^{\frac{1-2\gamma }{1-\gamma}} dx <\infty\,.
\end{equation}
Let $X_i=f(Y_i) - \nu_{\gamma}(f)$ where $f$ belongs to $\mathcal{F}^*(H,
\nu_{\gamma})$. Then (\ref{cond}) is satisfied and the conclusion of Theorem \ref{mainth}  holds with
\begin{equation}  \label{defeta*}
\eta= \nu_\gamma((f-\nu_\gamma(f))^2)+ 2 \sum_{k>0} \nu_\gamma
((f-\nu_\gamma(f))f\circ T_\gamma^k) \, .
\end{equation}
\end{cor}

\noindent \textbf{Proof.} To prove this corollary, it suffices to see that (\ref%
{lilcond}) implies (\ref{condalpha}). For
 this purpose, we use Proposition 1.17 in Dedecker {\it et al.} (2010) stating that there exist two positive constant $B,C$
such that, for any $n>0$, $Bn^{(\gamma -1)/\gamma }\leq \alpha _{\mathbf{Y}%
}(n)\leq Cn^{(\gamma -1)/\gamma }$, together with their computations page
817. $\square $

\medskip

In particular, if $f$ is $\BV$ and $\gamma<1/2$, we infer from Corollary \ref{ASmapB} that the conclusion of Theorem \ref{mainth} holds  with $\eta$ defined by (\ref{defeta*}) . Note also
that \eqref{lilcond} is satisfied if $H$ is such that $H(x)\leq C
x^{-2(1-\gamma)/(1-2\gamma)}(\ln(x))^{-b}$ for $x$ large enough and $%
b>(1-\gamma)/(1-2\gamma)$. Therefore, since the density $h_{\nu_{\gamma}}$
of $\nu_{\gamma}$ is such that $h_{\nu_{\gamma}}(x) \leq C x^{-\gamma}$ on $%
(0, 1]$, one can easily prove that if $f$ is positive and non increasing on
(0, 1), with $$f(x) \leq \frac{C}{x^{(1-2\gamma)/2}|\ln(x)|^{d}} \quad \text{near 0 for some
$d>1/2$},$$ then \eqref{lilcond} and the quenched FCLT hold. Notice that when $f$ is exactly of the
form $f(x)=x^{-(1-2\gamma)/2}$, Gou\"{e}zel (2004) proved that the central
limit theorem holds for $\sum_{i=1}^n ( f(Y_i) - \nu_{\gamma}(f)) $ but with
the normalization $\sqrt{n \ln (n)}$. This shows that the condition (\ref%
{lilcond}) is essentially optimal for the quenched CLT with the normalization $\sqrt n$.

\section{Some general results}

\quad In this section we develop sufficient conditions imposed to
conditional expectations of
partial sums  for the validity of the quenched CLT and FCLT.

For any positive integers $i$ and $p$,
define $S_{p}^{(i)}=S_{pi}-S_{p(i-1)}$.

\subsection{A quenched CLT}

Let us introduce the following three
conditions under which the quenched central limit theorem holds:
\begin{eqnarray*}
\text{\textbf{C$_{1}$}} &&\lim_{m\rightarrow \infty }\limsup_{p\rightarrow
\infty }\frac{1}{\sqrt{mp}}\sum_{i=2}^{m+1}\mathbb{E}_{0}|\mathbb{E}%
_{(i-2)p}(S_{p}^{(i)})|=0\ \ a.s. \\
\text{\textbf{C$_{2}$}} &&\text{there exists a $T$-invariant r.v. $\eta$ that is ${\mathcal{F}}_{0}$-measurable and such that} \\
&&\lim_{m\rightarrow \infty }\limsup_{p\rightarrow
\infty }\mathbb{E}_{0}\Big |\sum_{i=1}^{m}\frac{1}{mp}\mathbb{E}_{(i-1)p}%
\big ((S_{p}^{(i+1)})^{2}\big )-\eta \Big |=0\ \ a.s. \\
&&\lim_{m\rightarrow \infty }\limsup_{p\rightarrow \infty }\mathbb{E}_{0}%
\Big |\sum_{i=1}^{m}\frac{1}{mp}\mathbb{E}_{(i-1)p}\big (%
(S_{p}^{(i)}+S_{p}^{(i+1)})^{2}\big )-2\eta \Big |=0\ \ a.s. \\
\text{\textbf{C$_{3}$}} &&\text{for each $\varepsilon >0$}\quad
\lim_{m\rightarrow \infty }\limsup_{p\rightarrow \infty }\frac{1}{m}%
\sum_{i=1}^{m}\frac{1}{p}\mathbb{E}_{0}\big ((S_{p}^{(i)})^{2}\mathbf{1}%
_{|S_{p}^{(i)}|/\sqrt{p}>\varepsilon \sqrt{m}}\big )=0\ \ a.s.
\end{eqnarray*}

\begin{prop}
\label{mainprop} Assume that \textbf{C$_{1}$}, \textbf{C$_{2}$} and \textbf{C%
$_{3}$} hold. Then, on a set of probability one, for any continuous and bounded function $f$,
\begin{equation*}
\lim_{n\rightarrow \infty }{\mathbb{E}}_{0}(f(n^{-1/2}S_{n}))=\int f(x\sqrt{%
\eta })g(x)dx \, ,
\end{equation*}%
where $g$ is the density of a standard normal.
\end{prop}

This proposition is designed especially for the proof of Theorem \ref{mainth}.
 Notice that in the expression $\mathbb{E}_{(i-2)p}(S_{p}^{(i)})$ of
condition \textbf{C$_{1}$} there is a gap of $p$ variables between $%
S_{p}^{(i)}$ and the variables used for conditioning. This gap is important
for weakening the dependence and is essentially used in the proof of Theorem %
\ref{mainth}. 

\medskip

\noindent \textit{Proof of Proposition \ref{mainprop}.} The result will
follow from Proposition \ref{propgap2} below, for double indexed arrays of random
variables:

\begin{prop}
\label{propgap2}Assume that $(Y_{n,m,i})_{i\geq 1}$ is an array of random
variables in ${\mathbb{L}}^{2}$ adapted to an array
$(\mathcal{G}_{n,m,i})_{i\geq 1}$ of nested sigma fields. Let
${\mathbb E}_{n,m,i}$ denote the conditional expectation with respect
to $\mathcal{G}_{n,m,i}$.
Suppose that
\begin{equation}
\lim_{m\rightarrow \infty }\limsup_{n\rightarrow \infty }\sum_{i=2}^{m+1}%
\mathbb{E}|\mathbb{E}_{n,m,i-2}(Y_{n,m,i})|=0\,,  \label{C1}
\end{equation}%
and that there exists $\sigma ^{2}\geq 0$ such that
\begin{equation}
\lim_{m\rightarrow \infty }\limsup_{n\rightarrow \infty }\mathbb{E}%
\Big |\sum_{i=1}^{m}\mathbb{E}_{n,m,i-1}\big (Y_{n,m,i+1}^{2}\big )-\sigma ^{2}
\Big |=0%
\text{ }  \label{C2one}
\end{equation}%
and
\begin{equation}
\lim_{m\rightarrow \infty }\limsup_{n\rightarrow \infty }\mathbb{E}%
\Big|\sum_{i=1}^{m}\mathbb{E}_{n,m,i-1}\big ((Y_{n,m,i}+Y_{n,m,i+1})^{2}\big )%
-2\sigma ^{2}\Big|=0\,.  \label{C2bis}
\end{equation}%
Assume in addition that for each $\varepsilon >0$
\begin{equation}
\lim_{m\rightarrow \infty }\limsup_{n\rightarrow \infty }\ \sum_{i=1}^{m+1}%
\mathbb{E}(Y_{n,m,i}^{2}\mathbf{1}_{|Y_{n,m,i}|>\varepsilon })=0.  \label{L}
\end{equation}%
Then for any continuous and bounded function $f$,
\begin{equation*}
\lim_{m\rightarrow \infty }\limsup_{n\rightarrow \infty }\Big |\mathbb{E}%
\Big(f\Big(\sum_{i=1}^{m}Y_{n,m,i}\Big)\Big)-\mathbb{E}(f(\sigma N))\Big |=0\,,
\end{equation*}%
where $N$ is a standard Gaussian random variable.
\end{prop}

Before proving Proposition \ref{propgap2}, let us show how it leads to
Proposition \ref{mainprop}. Let $m$ be a fixed positive integer less than $n$.
Set $p=[n/m]$. We apply Proposition \ref{propgap2} to the sequence $%
Y_{n,m,i}=S_{p}^{(i)}/\sqrt{mp}$ and the filtration
${\mathcal{G}}_{n,m,i}={\mathcal{F}}_{ip}$. We also replace the expectation
${\mathbb E}$ by the conditional expectation ${\mathbb E}_0$ (recall that all the conditional expectations
of functions of $T$
with respect to ${\mathcal F}_0$ are obtained through the regular conditional probability
$P_{T|{\mathcal F}_0})$, and
$\sigma^2$ by the non negative ${\mathcal F}_0$-measurable random variable $\eta$. With these
notations,  the conditions \textbf{C$_{1}$}, \textbf{C$_{2}$} and \textbf{C$_{3}$}
imply that (\ref{C1}), (\ref{C2one}), (\ref{C2bis}) and $(\ref{L})$ hold  almost surely.
It follows from Proposition \ref{propgap2} that, on a set of probability one, for any continuous and bounded function $f$,
\begin{equation*}
\lim_{m\rightarrow \infty }\limsup_{n\rightarrow \infty }\Big |{\mathbb{E}}%
_{0}\Big (f\Big(n^{-1/2}\sum_{i=1}^{m[n/m]}X_{i}\Big)\Big )-\int f(x\sqrt{\eta }%
)g(x)dx\Big |=0\ \, ,
\end{equation*}%
where $g$ is the density of a standard normal.
Proposition \ref{mainprop} will then follow if we can prove that for any $%
\varepsilon >0$,
\begin{equation}
\lim_{m\rightarrow \infty }\limsup_{n\rightarrow \infty }{\mathbb{P}}_{0}%
\Big (\Big|\sum_{i=1}^{n}X_{i}-\sum_{i=1}^{m[n/m]}X_{i}\Big|\geq \varepsilon \sqrt{n}%
\Big)=0\ \text{\ }a.s.  \label{neglp0}
\end{equation}%
With this aim, we notice that
\begin{equation*}
{\mathbb{P}}_{0}\Big (\Big|\sum_{i=1}^{n}X_{i}-\sum_{i=1}^{m[n/m]}X_{i}\Big|\geq
\varepsilon \sqrt{n}\Big)\leq {\mathbb{P}}_{0}\Big (m^{2}\max_{1\leq i\leq
n}X_{i}^{2}\geq \varepsilon ^{2}n\Big).
\end{equation*}%
and therefore (\ref{neglp0}) holds by relation \eqref{maxERG} in Lemma \ref{ergodic} applied to $
Z_{i}=X_{i}^{2}$. It remains to prove Proposition \ref{propgap2}.

\medskip

\noindent \textit{Proof of Proposition \ref{propgap2}.}
For any positive integer $i$, let
\begin{equation}
U_{n,m,i}=Y_{n,m,i}+\mathbb{E}_{n,m,i}(Y_{n,m,i+1})-\mathbb{E}%
_{n,m,i-1}(Y_{n,m,i})\,.  \label{defxnmi}
\end{equation}%
To ease the notation, we shall drop the first two indexes (the pair $n,m)$
when no confusion is possible. With this notation,
\begin{equation*}
Y_{i}=U_{i}-\mathbb{E}_{i}(Y_{i+1})+\mathbb{E}_{i-1}(Y_{i})\,,
\end{equation*}%
and since we have telescoping sum,
\begin{equation*}
\sum_{i=1}^{m}Y_{i}=\sum_{i=1}^{m}U_{i}+\mathbb{E}_{0}(Y_{1})-\mathbb{E}%
_{m}(Y_{m+1})\,.
\end{equation*}%
Notice that for any $i\in \{1,m+1\}$ and any $\varepsilon >0$,
\begin{equation}
\mathbb{E}(|\mathbb{E}_{i-1}(Y_{i})|^{2})\leq \varepsilon ^{2}+\mathbb{E}%
(Y_{i}^{2}\mathbf{1}_{|Y_{i}|>\varepsilon })\,.
\end{equation}%
Therefore by condition (\ref{L}),
\begin{equation}
\lim_{m\rightarrow \infty }\limsup_{n\rightarrow \infty }\mathbb{E}\big (%
\mathbb{E}_{n,m,m}(Y_{n,m,m+1}))^{2}+(\mathbb{E}_{n,m.0}(Y_{n,m,1}))^{2}\big
)=0\,.  \label{p0gap2}
\end{equation}%
The theorem will be proven if we can show that the sequence $%
(U_{n,m,i})_{i\geq 1}$ defined by (\ref{defxnmi}) satisfies the conditions
of Theorem \ref{thmgap1}. We first notice that
$
\mathbb{E}_{i-1}(U_{i})=\mathbb{E}_{i-1}(Y_{i+1})
$.
Hence condition (\ref{C1ap}) is clearly satisfied under (\ref{C1}). On an
other hand,
\begin{multline}
\mathrm{Var}(U_{i}|\mathcal{G}_{i-1})=\mathbb{E}_{i-1}\big (Y_{i}^{2}+2Y_{i}
\mathbb{E}_{i}(Y_{i+1})\big )  \label{varcond}
+\mathbb{E}_{i-1}\big ((\mathbb{E}_{i}(Y_{i+1}))^{2}\big )-(\mathbb{E}
_{i-1}(Y_{i}))^{2}   \\
-2(\mathbb{E}_{i-1}(Y_{i}))(\mathbb{E}_{i-1}(Y_{i+1}))-(\mathbb{E}
_{i-1}(Y_{i+1}))^{2}\,.
\end{multline}%
Notice that for any $\varepsilon >0$
\begin{eqnarray}
 \sum_{i=1}^{m}\mathbb{E}\big ((\mathbb{E}_{i-1}(Y_{i+1}))^{2}\big )  \nonumber
& \leq & \varepsilon \sum_{i=1}^{m}\mathbb{E}\big |\mathbb{E}_{i-1}(Y_{i+1})
\big |+\varepsilon \sum_{i=1}^{m}\mathbb{E}\big (|Y_{i+1}|\mathbf{1}
_{|Y_{i+1}|>\varepsilon }\big )
 +\sum_{i=1}^{m}\mathbb{E}\big (Y_{i+1}^{2}\mathbf{1}_{|Y_{i+1}|>
\varepsilon }\big )  \nonumber \\
& \leq &\varepsilon \sum_{i=1}^{m}\mathbb{E}\big |\mathbb{E}_{i-1}(Y_{i+1})
\big |+2 \sum_{i=2}^{m+1}\mathbb{E}\big (Y_{i}^{2}\mathbf{1}
_{|Y_{i}|>\varepsilon }\big )\,.  \label{p1gap2}
\end{eqnarray}
Similarly, for any $\varepsilon >0$,
\begin{equation}
\sum_{i=1}^{m} \mathbb{E}\big |(\mathbb{E}_{i-1}(Y_{i}))(\mathbb{E}
_{i-1}(Y_{i+1}))\big |  \notag  \label{p2gap2}
 \leq \varepsilon \sum_{i=1}^{m}\mathbb{E}\big |\mathbb{E}_{i-1}(Y_{i+1})
\big |+2 \sum_{i=1}^{m+1}\mathbb{E}\big (Y_{i}^{2}\mathbf{1}
_{|Y_{i}|>\varepsilon }\big )\,.
\end{equation}%
In addition since
$
\mathbb{E}_{i-1}\big (Y_{i}^{2}+2Y_{i}\mathbb{E}_{i}(Y_{i+1})\big )=
\mathbb{E}_{i-1}\big ((Y_{i}+Y_{i+1})^{2}\big )-\mathbb{E}_{i-1}\big (
Y_{i+1}^{2}\big )
$,
the conditions (\ref{C2one}) and (\ref{C2bis}) imply that
\begin{equation}
\lim_{m\rightarrow \infty }\limsup_{n\rightarrow \infty }\mathbb{E}\Big |%
\sum_{i=1}^{m+1}\mathbb{E}_{n,m,i-1}\big (Y_{n,m,i}^{2}+2Y_{n,m,i}\mathbb{E}%
_{n,m,i}(Y_{n,m,i+1})\big )-\sigma ^{2}\Big |=0\,.  \label{C2}
\end{equation}%
Starting from (\ref{varcond}) and considering (\ref{p1gap2}), (\ref{p2gap2})
and (\ref{C2}), it follows that condition (\ref{C2ap}) will be satisfied
provided that (\ref{C1}) and (\ref{L}) hold and
\begin{equation}
\lim_{m\rightarrow \infty }\limsup_{n\rightarrow \infty }\mathbb{E}\Big |%
\sum_{i=1}^{m}\big (\mathbb{E}_{n,m,i-1}\big ((\mathbb{E}%
_{n,m,i}(Y_{n,m,i+1}))^{2}\big )-(\mathbb{E}_{n,m,i-1}(Y_{n,m,i}))^{2}\Big )%
\big |=0\,.  \label{p3gap2}
\end{equation}%
To prove (\ref{p3gap2}), we first write that
\begin{multline*}
\sum_{i=1}^{m} \big (\mathbb{E}_{i-1}\big ((\mathbb{E}_{i}(Y_{i+1}))^{2}
\big )-(\mathbb{E}_{i-1}(Y_{i}))^{2}\big )
 =\mathbb{E}_{m}(Y_{m+1}))^{2}-(\mathbb{E}_{0}(Y_{1}))^{2}
 \\ -\sum_{i=1}^{m}\big ((\mathbb{E}_{i}(Y_{i+1}))^{2}
 -\mathbb{E}_{i-1}\big ((
\mathbb{E}_{i}(Y_{i+1}))^{2}\big )\big )\,.
\end{multline*}
By (\ref{p0gap2}), it follows that (\ref{p3gap2}) will hold if we can show
that
\begin{equation}
\lim_{m\rightarrow \infty }\limsup_{n\rightarrow \infty }\mathbb{E}\Big |%
\sum_{i=1}^{m}\big ((\mathbb{E}_{n,m,i}(Y_{n,m,i+1}))^{2}-\mathbb{E}%
_{n,m,i-1}\big ((\mathbb{E}_{n,m,i}(Y_{n,m,i+1}))^{2}\big )\big )\Big |=0\,.
\label{p3gap2bis}
\end{equation}%
This follows from an application of Lemma \ref{LGNmart} with
\begin{equation*}
d_{n,m,i}=(\mathbb{E}_{n,m,i}(Y_{n,m,i+1}))^{2}-\mathbb{E}_{n,m,i-1}\big ((%
\mathbb{E}_{n,m,i}(Y_{n,m,i+1}))^{2}\big )\,.
\end{equation*}%
Indeed
\begin{equation*}
\sum_{i=1}^{m}\mathbb{E}(|d_{n,m,i}|)\leq 2\sum_{i=1}^{m+1}\mathbb{E}%
(Y_{n,m,i}^{2})\,,
\end{equation*}%
and by Lemma \ref{curiouslma}, for any $\varepsilon >0$,
\begin{multline*}
\sum_{i=1}^{m}\mathbb{E}(|d_{n,m,i}|\mathbf{1}_{|d_{n,m,i}|>8\varepsilon
^{2}})\leq 2\sum_{i=1}^{m}\mathbb{E}\big ((\mathbb{E}%
_{n,m,i}(Y_{n,m,i+1}))^{2}\mathbf{1}_{(\mathbb{E}%
_{n,m,i}(Y_{n,m,i+1}))^{2}>4\varepsilon ^{2}}) \\
\leq 2\sum_{i=1}^{m}\mathbb{E}\big (Y_{n,m,i+1}^{2}\mathbf{1}_{|\mathbb{E}%
_{n,m,i}(Y_{n,m,i+1})|>2\varepsilon })\leq 4\sum_{i=1}^{m+1}\mathbb{E}\big (%
Y_{n,m,i}^{2}\mathbf{1}_{|Y_{n,m,i}|>\varepsilon })\,.
\end{multline*}%
So condition (\ref{1}) holds by using (\ref{C2one}) and (\ref{L}).

It remains to prove that (\ref{Lap}) holds. Clearly this can be achieved by
using (\ref{L}) combined with Lemma \ref{curiouslma}. $\square$

\subsection{Finite dimensional convergence}

\quad For $0<t_{1}<\cdots <t_{d}\leq 1$, define the function $\pi
_{t_{1},\ldots ,t_{d}}$ from $C([0,1])$ to ${\mathbb{R}}^{d}$ by $\pi
_{t_{1},\ldots ,t_{d}}(x)=(x(t_{1}),x(t_{2})-x(t_{1}),\ldots
,x(t_{d})-x(t_{d-1}))$. For any $a$ in ${\mathbb{R}}^{d}$ define the
function $f_{a}$ from ${\mathbb{R}}^{d}$ to ${\mathbb{R}}$ by $%
f_{a}(x)=<a,x>=\sum_{i=1}^{d}a_{i}x_{i}$.

\begin{prop}
\label{fidi} Assume that \textbf{C$_{1}$}, \textbf{C$%
_{2}$} and \textbf{C$_{3}$} hold. Then, on a set of probability one, for any continuous and bounded
function $h$, for any $a \in {\mathbb Q}^d$ and any $t_{1},t_{2},\ldots ,t_{d}$  rational numbers such
that $0<t_{1}<\cdots <t_{d}\leq 1$,
\begin{equation} \label{resfiniconv}
\lim_{n\rightarrow \infty }{\mathbb{E}}_{0}\Big (h\circ f_{a}\circ
\pi_{t_{1},\ldots ,t_{d}}(W_{n})\Big )=\int h\circ f_{a}\circ
\pi_{t_{1},\ldots ,t_{d}}(z\sqrt{\eta })W(dz) \, ,
\end{equation}
where $W$ is the distribution of a standard Wiener process.
\end{prop}

\noindent \textit{Proof of Proposition \ref{fidi}.} Since $\bigcup_{ d =1}^{\infty} {\mathbb Q}^d$ is countable, it suffices to prove that for any $a \in {\mathbb R}^d$ and any $t_{1},t_{2},\ldots ,t_{d}$  rational numbers such
that $0<t_{1}<\cdots <t_{d}\leq 1$, on a set of probability one, for any continuous and bounded
function $h$, the convergence \eqref{resfiniconv} holds. With this aim, for any $\ell \in
\{1,\ldots ,d\}$, we set $t_{\ell }=r_{\ell }/s_{\ell }$ where $r_{\ell }$ and $%
s_{\ell }$ are positive integers. Let $c_{d}=\prod_{\ell =1}^{d}s_{\ell }$.
Rewrite $t_{\ell }=b_{\ell }/c_{d }$. The $b_{\ell }$'s are then positive
integers such that $0<b_{1}<\cdots <b_{d}\leq c_{d}$. Let $m$ be a fixed
positive integer and let $p=[n/(mc_{d})]$. Notice that for any $\ell \in
\{1,\ldots ,d\}$,
\begin{equation*}
\lbrack nt_{\ell }]-mb_{\ell }<mpb_{\ell }\leq \lbrack nt_{\ell }]+1\,.
\end{equation*}%
Therefore for any reals $a_{1},\cdots ,a_{d}$, with the convention that $%
t_{0}=0$ and $b_{0}=0$,
\begin{equation*}
\Big |\sum_{\ell =1}^{d}a_{\ell }\sum_{i=[nt_{\ell -1}]+1}^{[nt_{\ell
}]}X_{i}-\sum_{\ell =1}^{d}a_{\ell }\sum_{i=pmb_{\ell -1}+1}^{pmb_{\ell
}}X_{i}\Big |\leq \sum_{\ell =1}^{d}|a_{\ell }|\sum_{i=pmb_{\ell
}+1}^{(p+1)mb_{\ell }}|X_{i}|\,.
\end{equation*}%
Using (\ref{maxERG}) of Lemma \ref{ergodic}, we infer that for any $\ell \in
\{1,\ldots ,d\}$ and every $\varepsilon >0$,
\begin{equation*}
\lim_{n\rightarrow \infty }\mathbb{P}_{0}\Big (\frac{|a_{\ell }|}{\sqrt{n}}%
\sum_{i=pmb_{\ell }+1}^{(p+1)mb_{\ell }}|X_{i}|>\varepsilon \Big )=0\ \text{
}a.s.
\end{equation*}%
In addition,
\begin{equation*}
\Big |\sum_{\ell =1}^{d}a_{\ell }\big (W_{n}(t_{\ell })-W_{n}(t_{\ell -1})%
\big )-\sum_{\ell =1}^{d}a_{\ell }\big (S_{[nt_{\ell }]}-S_{[nt_{\ell -1}]}%
\big )\Big |\leq 2\sum_{\ell =1}^{d}|a_{\ell }|\max_{1\leq i\leq n}|X_{i}|\,,
\end{equation*}%
implying once again by (\ref{maxERG}) in Lemma \ref{ergodic} that
\begin{equation}
\lim_{n\rightarrow \infty }n^{-1/2}\mathbb{E}_{0}\Big (\Big |\sum_{\ell
=1}^{d}a_{\ell }\big (W_{n}(t_{\ell })-W_{n}(t_{\ell -1})\big )-\sum_{\ell
=1}^{d}a_{\ell }\big (S_{[nt_{\ell }]}-S_{[nt_{\ell -1}]}\big )\Big |\Big )%
=0\ a.s.  \label{convasmax}
\end{equation}%
From the preceding considerations, it remains to prove
 that, on a set of probability one, for any continuous and bounded function $f$,
\begin{equation}
\lim_{m\rightarrow \infty }\limsup_{n\rightarrow \infty }\Big |\mathbb{E}%
_{0}\Big(f\Big(n^{-1/2}\sum_{\ell =1}^{d}a_{\ell }\sum_{i=pmb_{\ell
-1}+1}^{pmb_{\ell }}X_{i}\Big)\Big)-\mathbb{E}_{0}(f(\sigma _{d}N))\Big |=0 \, ,
\label{res1proofpropfidi}
\end{equation}%
where $\sigma _{d}^{2}=\eta \sum_{\ell =1}^{d}a_{\ell }^{2}(t_{\ell
}-t_{\ell -1})$ and $N$ is a standard
Gaussian random variable independent of ${\mathcal F}_0$. With this aim, we
write that
\begin{equation*}
\sum_{\ell =1}^{d}a_{\ell }\sum_{i=pmb_{\ell -1}+1}^{pmb_{\ell
}}X_{i}=\sum_{\ell =1}^{d}a_{\ell }\sum_{i=mb_{\ell -1}+1}^{mb_{\ell
}}S_{p}^{(i)}=\sum_{k=1}^{mb_{d}}\lambda _{m,d,k}S_{p}^{(k)}\,,
\end{equation*}%
where $\lambda _{m,d,k}=\sum_{\ell =1}^{d}a_{\ell }\mathbf{1}_{mb_{\ell
-1}+1\leq k\leq mb_{\ell }}$. Hence to prove (\ref{res1proofpropfidi}), it
suffices to apply Proposition \ref{propgap2} to the random variables $%
Y_{n,m,i}=(mpc_{d})^{-1/2}\lambda _{m,d,i}S_{p}^{(i)}$ and the filtration $%
\mathcal{G}_{n,m,i}=\mathcal{F}_{ip}$, by replacing the expectation $\mathbb{%
E}$ by $\mathbb{E}_{0}$. The conditions (\ref{C1}) and (\ref{L}) are
verified by using respectively \textbf{C$_{1}$} and \textbf{C$_{3}$}. To
verify (\ref{C2one}) and (\ref{C2bis}) with $\sigma ^{2}=\sigma
_{d}^{2}=\eta \sum_{\ell =1}^{d}a_{\ell }^{2}(t_{\ell }-t_{\ell -1})$, we
proceed as follows. For (\ref{C2one}), we write that
\begin{eqnarray*}
\mathbb{E}_{0}\Big |\sum_{i=1}^{mb_{d}}\mathbb{E}_{n,m,i-1}(Y_{n,m,i+1}^{2})-
\sigma _{d}^{2}\Big |&=&\mathbb{E}_{0}\Big |\frac{1}{mpc_{d}}\sum_{\ell
=1}^{d}a_{\ell}^{2}\sum_{i=mb_{\ell -1}+1}^{mb_{\ell }}\mathbb{E}%
_{(i-1)p}((S_{p}^{(i+1)})^{2})-\sigma _{d}^{2}\Big | \\
&\leq & \sum_{\ell =1}^{d}a_{\ell}^{2}\mathbb{E}_{0}\Big |\frac{1}{mpc_{d}}
\sum_{i=mb_{\ell -1}+1}^{mb_{\ell }}\mathbb{E}%
_{(i-1)p}((S_{p}^{(i+1)})^{2})-\eta (t_{\ell }-t_{\ell -1})\Big | \, .
\end{eqnarray*}
Since $t_\ell=b_\ell/c_d$, we obtain that
\begin{eqnarray*}
\mathbb{E}_{0}\Big |\sum_{i=1}^{mb_{d}}\mathbb{E}_{n,m,i-1}(Y_{n,m,i+1}^{2})-
\sigma _{d}^{2}\Big |
&\leq & \sum_{\ell =1}^{d}\frac{a_{\ell}^{2}b_{\ell }}{c_{d}}\mathbb{E}_{0}\Big |%
\frac{1}{mpb_{\ell }}\sum_{i=1}^{mb_{\ell }}\mathbb{E}%
_{(i-1)p}((S_{p}^{(i+1)})^{2})-\eta \Big | \\
& & +\sum_{\ell =1}^{d}\frac{a_{\ell}^{2}b_{\ell -1}}{c_{d}}\mathbb{E}_{0}\Big |%
\frac{1}{mpb_{\ell -1}}\sum_{i=1}^{mb_{\ell -1}}\mathbb{E}%
_{(i-1)p}((S_{p}^{(i+1)})^{2})-\eta \Big |\,.
\end{eqnarray*}%
Condition (\ref{C2one}) is then proved by using the first part of \textbf{C$%
_{2}$}. Using similar arguments, we prove (\ref{C2bis}) by using the second
part of \textbf{C$_{2}$}. $\square $

\subsection{A  quenched invariance principle}

\quad Let us define the maximal version of \textbf{C$_{3}$}. For $k\leq l$,
let $\bar{S}_{k,l}=\max_{k\leq i\leq l}|S_{i}-S_{k}|$.
\begin{equation*}
\text{\textbf{C$_{4}$}}\quad \text{for any $\varepsilon >0$}\quad
\lim_{m\rightarrow \infty }\limsup_{p\rightarrow \infty }\frac{1}{m}%
\sum_{i=1}^{m}\frac{1}{p}\mathbb{E}_{0}\big (\bar{S}_{(i-1)p,ip}^{2}\mathbf{1%
}_{|\bar{S}_{(i-1)p,ip}|/\sqrt{p}>\varepsilon \sqrt{m}}\big )=0\ \ a.s.
\end{equation*}

\begin{prop}
\label{mainpropbis} Assume that \textbf{C$_1$}, \textbf{C$_2$} and \textbf{C$%
_4$} hold. Then, on a set of probability one, for any continuous and bounded function $f$ from $C([0,1])$
to ${\mathbb{R}}$,
\begin{equation*}
\lim_{n\rightarrow \infty} {\mathbb{E}}_0 (f(W_n) ) = \int f(x\sqrt{\eta}%
)W(dx) \, ,
\end{equation*}
where $W$ is the distribution of a standard Wiener process.
\end{prop}

\noindent \textit{Proof of Proposition \ref{mainpropbis}.} In this proof, $m$
will always denote a positive integer. Since \textbf{C$_{4}$} implies
\textbf{C$_{3}$}, it follows that Proposition \ref{fidi} holds. In what
follows, we shall prove that the process $\{W_{n}(t),t\in \lbrack 0,1]\}$ is
almost surely tight, that is, for any $\varepsilon >0$,
\begin{equation}
\lim_{m\rightarrow \infty }\limsup_{n\rightarrow \infty }{\mathbb{P}}_{0}%
\Big(\sup_{|t-s|\leq m^{-1}}|W_{n}(t)-W_{n}(s)|>\varepsilon \Big)=0\
\text{almost surely.}%
  \label{tightbillingsley}
\end{equation}
By standard arguments, (\ref{tightbillingsley}) together with Proposition \ref{fidi}
imply  Proposition \ref{mainpropbis}.

According to Inequality (25) in Brown (1971), to prove (\ref{tightbillingsley}) it suffices to show that, for any $\varepsilon >0$,
\begin{equation}
\lim_{m\rightarrow \infty }\limsup_{n\rightarrow \infty }\sum_{i=1}^{m}{%
\mathbb{P}}_{0}\Big(\sup_{(i-1)m^{-1}<t\leq
im^{-1}}|W_{n}(t)-W_{n}((i-1)m^{-1})|>\varepsilon \Big)=0\text{ \ }a.s.
\label{Brown}
\end{equation}%
Since $\sup_{t\in \lbrack
0,1]}|W_{n}(t)-n^{-1/2}S_{[nt]}|=n^{-1/2}\max_{1\leq i\leq n}|X_{i}|$, by
using (\ref{maxERG}) of Lemma \ref{ergodic}, it follows that (\ref{Brown})
is equivalent to
\begin{equation}
\lim_{m\rightarrow \infty }\limsup_{n\rightarrow \infty }\sum_{i=1}^{m}{%
\mathbb{P}}_{0}\Big(\sup_{(i-1)m^{-1}<t\leq
im^{-1}}|S_{[nt]}-S_{[n(i-1)m^{-1}]}|>\varepsilon \sqrt{n}\Big)=0\quad a.s.
\label{Brown2}
\end{equation}%
Let $p=[n/m]$, and note that, for any non negative integer $i$, $%
[nim^{-1}]-i<pi\leq \lbrack nim^{-1}]$. It follows that, for any integer $i$
in $[1,m]$,
$$
\sup_{(i-1)m^{-1}<t\leq im^{-1}}|S_{[nt]} -S_{[n(i-1)m^{-1}]}|
 \leq \bar{S}_{(i-1)p,ip}+\frac{1}{\sqrt{n}}
\sum_{k=[n(i-1)m^{-1}]-m}^{[n(i-1)m^{-1}]}|X_{k}|+\frac{1}{\sqrt{n}}
\sum_{k=[nim^{-1}]-m}^{[nim^{-1}]}|X_{k}|\,.
$$
Using (\ref{maxERG}) of Lemma \ref{ergodic}, we infer that
\begin{equation*}
\lim_{n\rightarrow \infty }\frac{1}{\sqrt{n}}%
\sum_{k=[n(i-1)m^{-1}]-m}^{[n(i-1)m^{-1}]}{\mathbb{E}}_{0}(|X_{k}|)=0\text{
\ }a.s.
\end{equation*}%
Hence, (\ref{Brown}) holds as soon as
\begin{equation*}
\lim_{m\rightarrow \infty }\limsup_{n\rightarrow \infty }\sum_{i=1}^{m}{%
\mathbb{P}}\Big(\bar{S}_{(i-1)p,ip}>\varepsilon \sqrt{n}\Big |{\mathcal{F}}%
_{0}\Big)=0\text{ \ }a.s.,
\end{equation*}%
which holds under \textbf{C$_{4}$}.

\section{Proof of Theorem \protect\ref{mainth} and additional comments}

\subsection{Proof of Theorem \protect\ref{mainth}}

We first prove that the series $\eta ={\mathbb{E}}(X_{0}^{2}|\mathcal{I})+2\sum_{k>0}{\mathbb{E}}(X_{0}X_{k}|%
\mathcal{I})  $ converges almost surely and in ${\mathbb{L}}^{1}$. With this aim, it suffices to prove that 
\begin{equation} \label{cvgeeta}
\sum_{k \geq 1} \Vert \E ( X_0 X_k | {\mathcal I} ) \Vert_1 < \infty \, .
\end{equation}
From Claim 1(b) in Dedecker and Rio (2000), $\E ( X_0 X_k | {\mathcal I} ) = \E ( \E ( X_0 X_k  | {\mathcal F}_{- \infty})| {\mathcal I} )$ almost surely, where ${\mathcal F}_{- \infty} = \bigcap_{k \in {\mathbb Z}} {\mathcal F}_{k}$. Hence 
$$
\Vert \E ( X_0 X_k | {\mathcal I} ) \Vert_1 \leq \Vert \E ( X_0 X_k | {\mathcal F}_{- \infty} ) \Vert_1 \leq \Vert X_0  \E_0 ( X_k  ) \Vert_1 \, ,
$$
which proves \eqref{cvgeeta} by using \eqref{cond}. 

We turn now to the rest of the proof.
 \begin{prop}
\label{mainresult} If (\ref{cond}) holds, then \textbf{C$_1$},\textbf{C$_2$}
and \textbf{C$_4$} hold, with $\eta$ defined in (\ref{eta}). In addition the
conclusion of Proposition \ref{mainpropbis} also holds for $f$ in ${\mathcal{%
H}}^*$.
\end{prop}

\noindent \textit{Proof of Proposition \ref{mainresult}.} We first prove
that the following reinforced version of \textbf{C$_{2}$} holds:  
\begin{eqnarray*}
\text{\textbf{C$_{2}^{\ast }$}} &&\text{there exists a $T$-invariant r.v. $\eta$ that is ${\mathcal{F}}_{0}$-measurable and such that} \\
&&\text{for any integer $i\geq 1$}\quad
\lim_{n\rightarrow \infty }\mathbb{E}_{0}\Big |\frac{1}{n}\mathbb{E}_{(i-1)n}%
\big ((S_{n}^{(i)})^{2}\big )-\eta \Big |=0\text{ \ }a.s. \text{ and } \\
&&\text{for any integer $i\geq 1$}\quad \lim_{n\rightarrow \infty }\mathbb{E}%
_{0}\Big |\frac{1}{n}\mathbb{E}_{(i-1)n}\big ((S_{n}^{(i)}+S_{n}^{(i+1)})^{2}%
\big )-2\eta \Big |=0\text{ \ }a.s.
\end{eqnarray*}

More precisely, we shall prove that \textbf{C$_{2}^{\ast }$} holds with $\eta$ defined in \eqref{eta}. We shall only prove the first part of \textbf{C$_2^*$}, the proof of the second
part being similar. For any positive integer $N$,
\begin{equation}  \label{start}
\frac{ \big ( S_{n}^{(i)} \big )^2}{ n}= \frac 1 n \sum_{j=(i-1)n
+1}^{in} X_j^2 +\frac 2 n \sum_{j=(i-1)n +1}^{in-1} \sum_{l=1}^{(in-j)\wedge
N}X_j X_{j+l} + R_{i,N} \, .
\end{equation}

Firstly,
\begin{equation*}
{\mathbb{E}}_{0}(|{\mathbb{E}}_{(i-1)n}(R_{i,N})|)\leq \frac{1}{n}%
\sum_{j=(i-1)n+1}^{in}{\mathbb{E}}_{0}\Big (\sum_{l>j+N}|X_{j}{\mathbb{E}}%
_{j}(X_{l})|\Big)\,.
\end{equation*}%
Let $Z_{j,N}=\sum_{l>j+N}|X_{j}{\mathbb{E}}_{j}(X_{l})|$ and note that, by
assumption, $Z_{j,N}=Z_{0,N}\circ T^{j}$ belongs to ${\mathbb{L}}^{1}$.
Applying the ergodic theorem in relation (\ref{ERG}) of Lemma \ref{ergodic}
we obtain that
\begin{equation*}
\lim_{n\rightarrow \infty }\frac{1}{n}\sum_{j=(i-1)n+1}^{in}{\mathbb{E}}%
_{0}(Z_{j,N})={\mathbb{E}}(Z_{0,N}|{\mathcal{I}})\text{ \ }a.s.
\end{equation*}%
Hence,
\begin{equation*}
\limsup_{n\rightarrow \infty }{\mathbb{E}}_{0}(|{\mathbb{E}}%
_{(i-1)n}(R_{i,N})|)\leq {\mathbb{E}}(Z_{0,N}|{\mathcal{I}})\ \ a.s.
\end{equation*}%
and consequently
\begin{equation}
\lim_{N\rightarrow \infty }\limsup_{n\rightarrow \infty }{\mathbb{E}}_{0}(|{%
\mathbb{E}}_{(i-1)n}(R_{i,N})|)=0\ \ a.s.  \label{residu}
\end{equation}

Next, let
\begin{eqnarray*}
\eta _{N} &=&{\mathbb{E}}(X_{0}^{2}|\mathcal{I})+2\sum_{k=1}^{N}{\mathbb{E}}%
(X_{0}X_{k}|\mathcal{I})\  \\
\text{and}\ \eta _{N,K} &=&{\mathbb{E}}(X_{0}^{2}\mathbf{1}_{|X_{0}|^{2}\leq
K}|\mathcal{I})+2\sum_{k=1}^{N}{\mathbb{E}}(X_{0}X_{k}\mathbf{1}%
_{|X_{0}X_{k}|\leq K}|\mathcal{I})\,.
\end{eqnarray*}%
By the ergodic theorem for stationary sequences,
\begin{equation}
\lim_{n\rightarrow \infty }\Big|\eta _{N,K}-\frac{1}{n}%
\sum_{j=(i-1)n+1}^{in}X_{j}^{2}\mathbf{1}_{|X_{j}|^{2}\leq K}-\frac{2}{n}%
\sum_{j=(i-1)n+1}^{in-1}\sum_{l=1}^{(in-j)\wedge N}X_{j}X_{j+l}\mathbf{1}%
_{|X_{j}X_{j+l}|\leq K}\Big |=0\text{ \ }a.s.  \label{erg1}
\end{equation}%
and by the ergodic theorem in relation (\ref{ERG}) of Lemma \ref{ergodic}
applied with $Z_{j}=X_{j}^{2}\mathbf{1}_{|X_{j}|^{2}>K}$ and with $%
Z_{j}=\sum_{l=1}^{N}|X_{j}X_{j+l}|\mathbf{1}_{|X_{j}X_{j+l}|>K}$,
\begin{equation}
\lim_{K\rightarrow \infty }\limsup_{n\rightarrow \infty }{\mathbb{E}}_{0}%
\Big(\frac{1}{n}\sum_{j=(i-1)n+1}^{in}X_{j}^{2}\mathbf{1}_{|X_{j}|^{2}>K}+
\frac{2}{n}\sum_{j=(i-1)n+1}^{in-1}\sum_{l=1}^{(in-j)\wedge N}|X_{j}X_{j+l}|%
\mathbf{1}_{|X_{j}X_{j+l}|>K}\Big)=0\text{ \ }a.s.  \label{erg2}
\end{equation}%
Using (\ref{erg1}), (\ref{erg2}) and the dominated convergence theorem, it
follows that
\begin{equation}
\lim_{n\rightarrow \infty }{\mathbb{E}}_{0}\Big(\Big |\eta _{N}-\frac{1}{n}%
\sum_{j=(i-1)n+1}^{in}X_{j}^{2}-\frac{2}{n}\sum_{j=(i-1)n+1}^{in-1}%
\sum_{l=1}^{(in-j)\wedge N}X_{j}X_{j+l}\Big |\Big)=0\ a.s.
\label{centralterm}
\end{equation}%
The first part of condition \textbf{C$_{2}^{\ast }$} follows from (\ref%
{start}), (\ref{residu}) and (\ref{centralterm}), and the fact that $%
\lim_{N\rightarrow \infty }\eta _{N}=\eta $ almost surely.

Next, we prove that \textbf{C$_{1}$} holds. With this aim, we first notice that it suffices to prove that for any
integer $i\geq 2$,
\begin{equation*}
\lim_{n\rightarrow \infty }{\mathbb{E}}(|X_{0}||{\mathcal{I}}){\mathbb{E}}%
_{0}\Big(\Big|{\mathbb{E}}_{(i-2)n}\Big(\frac{S_n^{(i)}}{\sqrt{n}}%
\Big)\Big|\Big)=0\ a.s.
\end{equation*}%
Indeed, on the invariant set where ${\mathbb{E}}(|X_{0}||{\mathcal{I}}) =0$ almost surely, the random variables $X_i$'s are equal to zero almost surely. 
Now, using the same arguments as in the proof of Lemma \ref{ergodic}, we can prove that ${\mathbb{E}}(|X_{0}||{
\mathcal{I}})=\E ( {\mathbb{E}}(|X_{0}||{
\mathcal{I}}) | {\mathcal F}_0)$ almost surely. Hence, for any integer $i \geq 2$, $${\mathbb{E}}(|X_{0}||{\mathcal{I}}){\mathbb{E}}%
_{0}\big(\big|{\mathbb{E}}_{(i-2)n}\big(S_n^{(i)}
\big)\big|\big) = {\mathbb{E}}_{0}\big(\big|{\mathbb{E}}_{(i-2)n}\big({\mathbb{E}}(|X_{0}||{
\mathcal{I}})S_n^{(i)}\big)\big|\big) \ a.s.$$
Now
\begin{multline}
\frac{1}{\sqrt{n}}{\mathbb{E}}_{0}\big(\big|{\mathbb{E}}_{(i-2)n}\big({\mathbb{E}}(|X_{0}||{
\mathcal{I}})S_n^{(i)}\big)\big|\big)
\label{merciflo} \\
\leq {\mathbb{E}}_{0}\Big(\Big|{\mathbb{E}}_{(i-2)n}\Big(\Big(\frac{1}{n}
\sum_{k=(i-2)n+1}^{(i-1)n}|X_{k}|-{\mathbb{E}}(|X_{0}||{\mathcal{I}})\Big)
\frac{S_n^{(i)}}{\sqrt{n}}\Big)\Big|\Big)   \\
+\frac{1}{n^{3/2}}\sum_{k=(i-2)n+1}^{(i-1)n}{\mathbb{E}}_{0}\big(\big |{\mathbb{E}}
_{(i-2)n}\big(|X_{k}| \, S_n^{(i)}\big)\big|\big)\,.
\end{multline}

Using the fact that $\mathcal{F}_{0}\subseteq {\mathcal{F}}_{(i-2)n}$ for any $i \geq 2$, and
applying Cauchy-Schwarz's inequality conditionally to ${\mathcal{F}}_{0}$,
the first term on right hand in (\ref{merciflo}) is smaller than
\begin{equation}
{\mathbb{E}}_{0}^{1/2}\Big(\Big(\frac{1}{n}\sum_{k=(i-2)n+1}^{(i-1)n}|X_{k}|-%
{\mathbb{E}}(|X_{0}||{\mathcal{I}})\Big)^{2}\Big){\mathbb{E}}_{0}^{1/2}\Big(%
\Big(\frac{S_n^{(i)}}{\sqrt{n}}\Big)^{2}\Big)\,.  \label{bound}
\end{equation}%
By \textbf{C$_{2}^{\ast }$},
\begin{equation}
\lim_{n\rightarrow \infty }{\mathbb{E}}_{0}\Big(\Big(\frac{S_n^{(i)}%
}{\sqrt{n}}\Big)^{2}\Big)=\eta \text{ \ }a.s.  \label{lim1}
\end{equation}%
Since $X_{0}$ belongs to ${\mathbb{L}}^{2}$, proceeding
as in the proof of (\ref{centralterm}), we obtain that
\begin{equation}
\lim_{n\rightarrow \infty }{\mathbb{E}}_{0}\Big(\Big(\frac{1}{n}%
\sum_{k=(i-2)n+1}^{(i-1)n}|X_{k}|-{\mathbb{E}}(|X_{0}||{\mathcal{I}})\Big)%
^{2}\Big)=0\ \text{a.s.}  \label{lim2}
\end{equation}%
From (\ref{bound}), (\ref{lim1}) and (\ref{lim2}), we infer that the first
term on right hand in (\ref{merciflo}) converges to $0$ almost surely as $n$
tends to infinity.

Now, for any integer $k$ belonging to $](i-2)n,(i-1)n]$,
\begin{align*}
\frac{1}{\sqrt{n}}{\mathbb{E}}_0\big(\big |{\mathbb{E}} \big(|X_{k}| \, S_n^{(i)}
\big|{\mathcal{F}}_{(i-2)n}\big)\big|\big) & \leq \frac{1}{\sqrt{n}}%
{\mathbb{E}}_0\big(\big |{\mathbb{E}}(|X_{k}|S_n^{(i)}|{\mathcal{F}%
}_{k})\big| \big) \\
&\leq \frac{1}{\sqrt{n}}{\mathbb{E}}_{0}\Big(\sum_{i=k+1}^{\infty }|X_{k}{%
\mathbb{E}}_{k}(X_{i})|\Big)\,.
\end{align*}%
Let $Z_{k}=\sum_{i=k+1}^{\infty }|X_{k}{\mathbb{E}}_{k}(X_{i})|$ and note
that, by assumption, $Z_{k}=Z_{0}\circ T^{k}$ belongs to ${\mathbb{L}}^{1}$.
It follows that the second term on the right-hand side of (\ref{merciflo}) is
smaller than $n^{-3/2}\sum_{k=(i-2)n+1}^{(i-1)n}{\mathbb{E}}_{0}(Z_{k})\,,$
which converges almost surely to $0$ as $n$ tends to infinity, by the
ergodic theorem in relation (\ref{ERG}) of Lemma \ref{ergodic}. Hence
\textbf{C$_{1}$} is proved.

We turn now to the proof of \textbf{C$_4$}. With this aim, we shall prove
the following reinforcement of it:
\begin{equation*}
\text{\textbf{C$_4^*$}} \quad \quad \quad \lim_{k \rightarrow \infty}
\limsup_{n \rightarrow \infty} \max_{1\leq i \leq k} {\mathbb{E}}_0\Big(%
\frac{\bar S_{(i-1)n, in}^2}{n} \Big(1 \wedge \frac{\bar S_{(i-1)n, in}}{%
\sqrt{nk}}\Big) \Big)= 0 \ a. s.
\end{equation*}
To prove \textbf{C$_4^*$}, we shall use the following maximal inequality,
which is a conditional version of the inequality given in Proposition 1(a) of Dedecker and Rio (2000).

\begin{prop}
\label{maxineq}For any $k<l$ and $\lambda \geq 0$
let  $
\Gamma_{k,l}(\lambda)= \{\bar S_{k,l}> \lambda\}$. The following inequality
holds
\begin{equation*}
{\mathbb{E}}_{0}((\bar{S}_{k,l}-\lambda )_{+}^{2})\leq 8\sum_{i=k+1}^{l}{
\mathbb{E}}_{0}(X_{i}^{2}\mathbf{1}_{\Gamma _{k,i}(\lambda
)})+16\sum_{i=k+1}^{l}{\mathbb{E}}_{0}(|X_{i}\mathbf{1}_{\Gamma
_{k,i}(\lambda )}{\mathbb{E}}_{i}(S_{l}-S_{i})|)\,.
\end{equation*}
\end{prop}

Let us continue the proof of \textbf{C$_{4}^{\ast }$}. Note first that
\begin{equation*}
{\mathbb{E}}_{0}\Big(\frac{\bar{S}_{(i-1)n,in}^{2}}{n}\Big(1\wedge \frac{%
\bar{S}_{(i-1)n,in}}{\sqrt{nk}}\Big)\Big)\leq 2\varepsilon {\mathbb{E}}_{0}%
\Big(\frac{\bar{S}_{(i-1)n,in}^{2}}{n}\Big)+\frac{4}{n}{\mathbb{E}}_{0}\Big((%
\bar{S}_{(i-1)n,in}-\varepsilon \sqrt{nk})_{+}^{2}\Big)\,.
\end{equation*}%
From Proposition \ref{maxineq} with $\lambda =0$, we obtain that
\begin{equation*}
{\mathbb{E}}_{0}\Big(\frac{\bar{S}_{(i-1)n,in}^{2}}{n}\Big)\leq \frac{8}{n}%
\sum_{k=(i-1)n+1}^{in}{\mathbb{E}}_{0}(X_{k}^{2})+\frac{16}{n}%
\sum_{k=(i-1)n+1}^{in}{\mathbb{E}}_{0}(Z_{k})\,,
\end{equation*}%
and, by the ergodic theorem in relation (\ref{ERG}) of Lemma \ref{ergodic},
\begin{equation}
\limsup_{n\rightarrow \infty }\max_{1\leq i\leq k}{\mathbb{E}}_{0}\Big(\frac{%
\bar{S}_{(i-1)n,in}^{2}}{n}\Big)\leq 8{\mathbb{E}}(X_{0}^{2}|{\mathcal{I}}%
)+16{\mathbb{E}}(Z_{0}|{\mathcal{I}})\ a.s.  \label{bounded}
\end{equation}%
Hence \textbf{C$_{4}^{\ast }$} will be proved if, for any $\varepsilon >0$,
\begin{equation}
\lim_{k\rightarrow \infty }\limsup_{n\rightarrow \infty }\max_{1\leq i\leq k}%
\frac{1}{n}{\mathbb{E}}_{0}\Big((\bar{S}_{(i-1)n,in}-\varepsilon \sqrt{nk}%
)_{+}^{2}\Big)=0\ a.s.  \label{squifaut}
\end{equation}
Applying Proposition \ref{maxineq}, we infer that, for any positive
integer $N$,
\begin{gather}
{\mathbb{E}}_{0}\Big((\bar{S}_{(i-1)n,in}-\varepsilon \sqrt{nk})_{+}^{2}\Big)%
\leq \frac{4}{n}\sum_{j=(i-1)n+1}^{in}{\mathbb{E}}_{0}(X_{j}^{2}\mathbf{1}%
_{\Gamma _{(i-1)n,in}(\varepsilon \sqrt{nk})})  \label{step1} \\
+\frac{8}{n}\sum_{j=(i-1)n+1}^{in}\sum_{l=1}^{(in-j)\wedge N}{\mathbb{E}}%
_{0}(|X_{j}X_{j+l}|\mathbf{1}_{\Gamma _{(i-1)n,in}(\varepsilon \sqrt{nk})})+%
\frac{8}{n}\sum_{j=(i-1)n+1}^{in}{\mathbb{E}}_{0}(Z_{j,N})\,,  \notag
\end{gather}
where $Z_{j,N}=\sum_{l>j+N}|X_{j}{\mathbb{E}}_{j}(X_{l})|$.  Since by \eqref{cond}, $Z_{j,N}=Z_{0,N} \circ T^j$ belongs to ${\mathbb{L}}^{1}$, The ergodic theorem in relation (\ref{ERG}) of Lemma \ref{ergodic} gives: 
for any positive integer $i$,
\begin{equation*}
\lim_{n\rightarrow \infty }\frac{1}{n}\sum_{j=(i-1)n+1}^{in}{\mathbb{E}}%
_{0}(Z_{j,N})={\mathbb{E}}(Z_{0,N}|{\mathcal{I}})\ a.s.
\end{equation*}%
and consequently,
\begin{equation}
\lim_{N\rightarrow \infty }\limsup_{n\rightarrow \infty }\max_{1\leq i\leq k}%
\frac{1}{n}\sum_{j=(i-1)n+1}^{in}{\mathbb{E}}_{0}(Z_{j,N})=0\ a.s.
\label{step2}
\end{equation}%
Now, for any positive $M$ and any $0\leq l\leq N$,
\begin{multline}
\frac{1}{n}\sum_{j=(i-1)n+1}^{in}{\mathbb{E}}_{0}(|X_{j}X_{j+l}|\mathbf{1}
_{\Gamma _{(i-1)n,in}(\varepsilon \sqrt{nk})})\leq \frac{M}{\varepsilon ^{2}k
}{\mathbb{E}}_{0}\Big(\frac{\bar{S}_{(i-1)n,in}^{2}}{n}\Big)  \label{step3}
\\
+\frac{1}{n}\sum_{j=(i-1)n+1}^{in}{\mathbb{E}}_{0}(|X_{j}X_{j+l}|\mathbf{1}
_{|X_{j}X_{j+l}|>M})\,.
\end{multline}
According to (\ref{bounded}), we have that
\begin{equation}
\lim_{k\rightarrow \infty }\limsup_{n\rightarrow \infty }\max_{1\leq i\leq k}%
\frac{M}{\varepsilon ^{2}k}{\mathbb{E}}_{0}\Big(\frac{\bar{S}_{(i-1)n,in}^{2}%
}{n}\Big)=0\ a.s.  \label{step4}
\end{equation}%
Next, by the ergodic theorem in relation (\ref{ERG}) of Lemma \ref{ergodic},
for any positive integer $i$,
\begin{equation*}
\lim_{n\rightarrow \infty }\frac{1}{n}\sum_{j=(i-1)n+1}^{in}{\mathbb{E}}%
_{0}(|X_{j}X_{j+l}|\mathbf{1}_{|X_{j}X_{j+l}|>M})={\mathbb{E}}(|X_{0}X_{l}|%
\mathbf{1}_{|X_{0}X_{l}|>M}|{\mathcal{I}})\ a.s.
\end{equation*}%
and consequently
\begin{equation}
\lim_{M\rightarrow \infty } \limsup_{k \rightarrow \infty}
\limsup_{n\rightarrow \infty }\max_{1\leq i\leq k}%
\frac{1}{n}\sum_{j=(i-1)n+1}^{in}{\mathbb{E}}_{0}(|X_{j}X_{j+l}|\mathbf{1}%
_{|X_{j}X_{j+l}|>M})=0\ a.s.  \label{step5}
\end{equation}%
Gathering (\ref{step1}), (\ref{step2}), (\ref{step3}), (\ref{step4}) and (%
\ref{step5}), we infer that (\ref{squifaut}) holds. This ends the proof of
\textbf{C$_{4}^{\ast }$}. Then, on a set of probability one, for any continuous and bounded function $\varphi $
from $C([0,1])$ to ${\mathbb{R}}$, \eqref{L1ps} follows by applying Proposition \ref{mainpropbis}. To prove that \eqref{L1ps}  also holds for $\varphi $ in ${\mathcal{H}}^{\ast }$, it suffices to notice that since (\ref{cond}) implies \textbf{C$%
_{4}^{\ast }$}, it entails in particular that almost surely, the sequence $%
(n^{-1}\max_{1\leq k\leq n}S_{k}^{2})_{n\geq 1}$ is uniformly integrable for
the conditional expectation with respect to ${\mathcal{F}}_{0}$. $\square $.

\medskip

\noindent \textit{Proof of Proposition \ref{maxineq}.} It is exactly the
same as to get (3.12) in the paper by Dedecker and Rio (2000), with the only
difference that the expectation is replaced by the conditional expectation
with respect to ${\mathcal{F}}_0$. $\square$.

\subsection{Some remarks on martingale approximations} \label{sectionremark}

The aim of this subsection is to point out that the conditions
\textbf{C$_{1}$}, \textbf{C$_{2}$} and \textbf{C$_{3}$} are
satisfied if there is an almost sure conditional martingale approximation in ${\mathbb L}^2$. This is another way to see that our conditions \textbf{C$_{1}$}, \textbf{C$_{2}$} and \textbf{C$_{3}$}
lead to sharp sufficient conditions for the quenched CLT.

From the proof of Theorem \ref{mainresult}, we see that, if $X_1$ is a martingale difference, that is ${\mathbb E}(X_1|{\mathcal F}_0)=0$ a.s.,
then the conditions \textbf{C$_{1}$}, \textbf{C$_{2}^{\ast }$} and \textbf{C$_{4}^{\ast }$} are satisfied. The following
claim is then easily deduced.  

\begin{claim} Let $X_0$ and $d_0$ be two ${\mathcal F}_0$-measurable, centered and square integrable random
variables with ${\mathbb E}(d_0\circ T|{\mathcal F}_0)=0$ a.s., and let $X_i=X_0 \circ T^i$ and $d_i=d_0\circ T^i$. Let $S_n=X_1+ \cdots + X_n$
and $M_n=d_1+ \cdots + d_n$.
\begin{enumerate}
\item If
$$
  \lim_{n \rightarrow \infty} \frac{1}{n}{\mathbb E}_0((S_n-M_n)^2)=0 \quad \text{almost surely,}
$$
then the conditions \textbf{C$_{1}$}, \textbf{C$_{2}$} and \textbf{C$_{3}$} are satisfied with $\eta=\E (d_0^2 |{\mathcal I})$.
\item
If
\begin{equation}\label{maxapprox}
  \lim_{n \rightarrow \infty} \frac{1}{n}{\mathbb E}_0\Big(\max_{1\leq k \leq n}(S_k-M_k)^2\Big)=0 \quad \text{almost surely,}
\end{equation}
then the conditions \textbf{C$_{1}$}, \textbf{C$_{2}$} and \textbf{C$_{4}$} are satisfied with $\eta=\E (d_0^2 |{\mathcal I})$.
\end{enumerate}
\end{claim}
In particular, if the condition of Maxwell and Woodroofe (2000) is satisfied
\begin{equation} \label{condMW*}
  \sum_{n>0}\frac{\|{\mathbb E}_0(S_n)\|_2}{n^{3/2}} < \infty \, ,
\end{equation}
then it follows from Cuny and Merlev\`ede (2012) that (\ref{maxapprox})
holds, so that the conditions \textbf{C$_{1}$}, \textbf{C$_{2}$} and \textbf{C$_{4}$}
are satisfied. We already know from Peligrad and Utev (2006) that the Maxwell and Woodroofe condition is sharp in some sense for the FCLT, and therefore for the quenched FCLT also. This shows that the conditions \textbf{C$_{1}$}, \textbf{C$_{2}$} and \textbf{C$_{4}$} are essentially sharp for the quenched FCLT.

We mention that the Maxwell and Woodroofe condition and our condition (\ref{cond}) are of independent interests. For instance, when applied to strongly mixing sequences the condition \eqref{condMW*} leads to sub-optimal results as pointed out in Merlev\`ede {\it et al.} (2006). Obviously, the same remark is true when we apply it to $\alpha$-dependent sequences as defined in Section \ref{sectionappli}. More precisely, this gives the condition: $\sum_{k\geq 0} (k+1)^{-1/2}\Big ( \int_{0}^{\alpha_{\mathbf{Y}} (k)}Q^{2}(u)du \Big )^{1/2}<\infty$  instead of \eqref{condalpha}. Hence,  when applied to non necessarily bounded functions of the Markov chain associated to the intermittent map given in Section 3.1, the criterion \eqref{condMW*} is satisfied as soon as $f$ belongs to $\mathcal{F}^*(H,\nu_{\gamma})$ and $H$ is such that $H(x)\leq C
x^{-2(1-\gamma)/(1-2\gamma)}(\ln(x))^{-b}$ for $x$ large enough and $b>2(1-\gamma)/(1-2\gamma)$. Recall that by  condition \eqref{condalpha}, we only need $b>(1-\gamma)/(1-2\gamma)$. In addition, Point (v) of the main theorem in Durieu and Voln\'y (2008) shows that one can find a stationary sequence $(X_i)_{ i \in {\mathbb Z}}$ adapted to an increasing and stationary filtration $({\mathcal F}_i)_{ i \in {\mathbb Z}}$ in such a way that the condition \eqref{condMW*} holds but $X_0 \E_0 (S_n)$ does not converge in ${\mathbb L}^1$ and so the condition \eqref{cond} fails. Analyzing the examples given in their paper, one can also prove that there are stationary sequences for which \eqref{cond} holds but \eqref{condMW*} does not. We can even say more: there are stationary sequences for which \eqref{cond} holds but not \eqref{condMW*}, neither the Gordin criterion \eqref{condGordinL1}, nor the Hannan-Heyde condition are satisfied. Recall that the Hannan-Heyde condition is the following: 
\begin{equation} \label{HHcond} 
\E(X_0|{\mathcal F}_{- \infty})=0 \ a.s. \ \text{ and } \ \sum_{n \geq 0} \Vert {\mathbb E}_0 (X_n)-{\mathbb E}_{-1} (X_n)\Vert_2 <  \infty \, ,
\end{equation}
where ${\mathcal F}_{- \infty} = \bigcap_{k \in {\mathbb Z}} {\mathcal F}_{k}$.

In what follows $(\Omega, {\mathcal A }, \mu)$ is a probability space and $T : \Omega \rightarrow \Omega$ a bijective bimeasurable transformation preserving the measure $\mu$. Then $(\Omega, {\mathcal A }, \mu , T)$ is called a dynamical system. We refer to Sina\u\i \,  (1962) for a precise definition of the entropy of a dynamical system, and for the properties of  dynamical systems with positive entropy. The proof  of the next proposition is given in the Appendix.
\begin{prop} \label{propDV} Let $(\Omega, {\mathcal A }, \mu , T)$ be an ergodic dynamical system with positive entropy. Let ${\mathcal F} \subset {\mathcal A}$ be a $T$-invariant $\sigma$-algebra, i.e. ${\mathcal F} \subset T^{-1} (\mathcal F)$ and let ${\mathcal F}_i = T^{-i} (\mathcal F)$. There exists a ${\mathcal F}_0$-measurable and centered function $f$ in ${\mathbb L}^2$ such that, setting $X_i = f \circ T^i$,  the condition \eqref{cond} is satisfied but the conditions \eqref{condGordinL1}, \eqref{condMW*} and \eqref{HHcond} fail.
\end{prop}
To be complete, note that  a stationary sequence can be constructed in such a way that \eqref{condGordinL1} holds but the condition \eqref{cond} fails (see Section 5.2 in Durieu and Voln\'y (2008)). Moreover, a stationary sequence can be constructed in such a way that the condition \eqref{HHcond} holds but the condition \eqref{cond} fails (see Theorem 1 in Durieu (2009)).

\section{Normal approximation  for double indexed arrays and auxiliary results}
\label{NA}

There are many situations when we are dealing with double indexed
sequences of random variables. For instance at each point in the two
dimensional space we start a random walk. Our motivation for this section
comes from the fact that in our blocking procedure we introduce a new
parameter, the number of blocks, $m$, that is kept fixed at the beginning.

The next theorem treats the martingale approximation for  double arrays
of random variables.

\begin{thm}
\label{thmgap1} Assume that $(U_{n,m,i})_{i\geq 1}$ is an array of random
variables in ${\mathbb{L}}^{2}$ adapted to an array $(\mathcal{G}%
_{n,m,i})_{i\geq 1}$ of nested sigma fields. Let
${\mathbb E}_{n,m,i}$ denote the conditional expectation with respect
to $\mathcal{G}_{n,m,i}$.
Suppose that
\begin{equation}
\lim_{m\rightarrow \infty }\limsup_{n\rightarrow \infty }\mathbb{E}\Big |%
\sum_{i=1}^{m}\mathbb{E}_{n,m,i-1}(U_{n,m,i})\Big |=0\,,  \label{C1ap}
\end{equation}%
there exists $\sigma ^{2}\geq 0$ such that
\begin{equation}
\lim_{m\rightarrow \infty }\limsup_{n\rightarrow \infty }\mathbb{E}\Big |%
\sum_{i=1}^{m}\mathrm{Var}(U_{n,m,i}|\mathcal{G}_{n,m,i-1})-\sigma ^{2}\Big |%
=0 \, , \label{C2ap}
\end{equation}%
and for each $\varepsilon >0$
\begin{equation}
\lim_{m\rightarrow \infty }\limsup_{n\rightarrow \infty }\ \sum_{i=1}^{m}%
\mathbb{E}(U_{n,m,i}^{2}\mathbf{1}_{|U_{n,m,i}|>\varepsilon })=0  \, . \label{Lap}
\end{equation}%
Then for any continuous and bounded function $f$,
\begin{equation}
\lim_{m\rightarrow \infty }\limsup_{n\rightarrow \infty }\Big |\mathbb{E}%
\Big(f\Big(\sum_{i=1}^{m}U_{n,m,i}\Big)\Big)-\mathbb{E}(f(\sigma N))\Big |=0\,,
\label{resthmgap1}
\end{equation}%
where $N$ is a standard Gaussian variable.
\end{thm}

\noindent \textbf{Proof of Theorem \ref{thmgap1}.} For any $i\geq 1$, let $d_{n,m,i}=U_{n,m,i}-%
\mathbb{E}_{n,m,i-1}(U_{n,m,i})$. By condition (\ref{C1ap}), the theorem
will follow if we can prove that (\ref{resthmgap1}) holds with $%
\sum_{i=1}^{m}d_{n,m,i}$ replacing $\sum_{i=1}^{m}U_{n,m,i}$. If $\sigma
^{2}=0$ the theorem is trivial. So we can assume without loss of generality
that $\sigma ^{2}=1$.
In the rest of the proof, in order to ease the notations,
we shall drop the first two indexes $(n,m)$, keeping them only when it is
necessary to avoid confusion. Let $\varepsilon $ and $M$ be positive reals
fixed for the moment. For any $i\geq 1$, let
\begin{equation*}
V_{i}=\sum_{\ell =1}^{i}\mathbb{E}_{\ell -1}(d_{\ell }^{2})\,
\ \text{ and } \
\,Y_{i}=d_{i}\mathbf{1}_{|d_{i}|\leq \varepsilon }\mathbf{1}_{V_{i}\leq M}\,.
\end{equation*}%
Notice first that
\begin{eqnarray*}
{\mathbb{P}} \Big (\sum_{i=1}^{m}d_{i}\neq \sum_{i=1}^{m}Y_{i}\Big )&\leq &{
\mathbb{P}}\Big(\max_{1\leq i\leq m}d_{i}>\varepsilon \Big)+{\mathbb{P}}(V_{m}>M) \\
&\leq & \frac{1}{\varepsilon ^{2}}\sum_{i=1}^{m}\mathbb{E}(d_{i}^{2}\mathbf{1}
_{|d_{i}|>\varepsilon })+\frac{1}{M}\Big (1+\mathbb{E}\Big |\sum_{i=1}^{m}
\mathrm{Var}(U_{i}|\mathcal{G}_{i-1})-1\Big |\Big )\,.
\end{eqnarray*}
Hence using Lemma \ref{curiouslma}, we get that    
$$
{\mathbb{P}} \Big (\sum_{i=1}^{m}d_{i}\neq \sum_{i=1}^{m}Y_{i}\Big )
 \leq \frac{12}{\varepsilon ^{2}}\sum_{i=1}^{m}\mathbb{E}(U_{i}^{2}\mathbf{1}
_{|U_{i}|>\varepsilon /4})+\frac{1}{M}\Big (1+\mathbb{E}\Big |\sum_{i=1}^{m}
\mathrm{Var}(U_{i}|\mathcal{G}_{i-1})-1\Big |\Big )\,.
$$
Therefore using (\ref{Lap}) and (\ref{C2ap}), it follows that for all $\varepsilon >0$, 
\begin{equation}
\limsup_{m\rightarrow \infty }\limsup_{n\rightarrow \infty }{\mathbb{P}}%
\Big (\sum_{i=1}^{m}d_{n,m,i}\neq \sum_{i=1}^{m}Y_{n,m,i}\Big ) \leq \frac{1}{M}\,.
\label{p1thmgap1}
\end{equation}%
We notice now that since $\mathbb{E}_{i-1}(d_{i})=0$ a.s. and $V_{i}$ is $%
\mathcal{G}_{i-1}$-measurable,
\begin{equation*}
\sum_{i=1}^{m}\mathbb{E}_{i-1}(Y_{i})=\sum_{i=1}^{m}\mathbf{1}_{V_{i}\leq M}%
\mathbb{E}_{i-1}(d_{i}\mathbf{1}_{|d_{i}|>\varepsilon })\,.
\end{equation*}%
Therefore by Lemma \ref{curiouslma},  
\begin{equation*}
\mathbb{E}\Big |\sum_{i=1}^{m}\mathbb{E}_{i-1}(Y_{i})\Big |\leq \frac{12}{%
\varepsilon }\sum_{i=1}^{m}\mathbb{E}(U_{i}^{2}\mathbf{1}_{|U_{i}|>%
\varepsilon /4})\,,
\end{equation*}%
implying, by using (\ref{Lap}), that for all positive reals $\varepsilon $ and $M$,
\begin{equation}
\lim_{m\rightarrow \infty
}\limsup_{n\rightarrow \infty }\mathbb{E}\Big |\sum_{i=1}^{m}\mathbb{E}%
_{i-1}(Y_{i})\Big |=0\,.  \label{p2thmgap1}
\end{equation}%
Considering (\ref{p1thmgap1}) and (\ref{p2thmgap1}), the theorem will follow
if we can show that for any continuous  bounded function $f$,
\begin{equation}
\lim_{M\rightarrow \infty }\limsup_{\varepsilon \rightarrow
0}\limsup_{m\rightarrow \infty }\limsup_{n\rightarrow \infty }\Big |\mathbb{E%
}\Big(f\Big(\sum_{i=1}^{m}d_{n,m,i}^{\ast  }\Big)\Big)-\mathbb{E}(f(N))\Big |=0\,,
\label{p3thmgap1}
\end{equation}%
where
\begin{equation*}
d_{n,m,i}^{\ast }=Y_{n,m,i}-\mathbb{E}_{n,m,i-1}(Y_{n,m,i})\,.
\end{equation*}%
Let 
$$
s^2_{n,m} =\sum_{i=1}^{m} \mathbb{E} ( d_{n,m,i}^{\ast 2} ) \, ,
$$
and notice that for any $\delta >0$, $\mathbb{E}\big(|d_{n,m,i}^{\ast }|^{2+2\delta }) < \infty$, $i=1,2,\dots$ Hence, by the first theorem stated in Heyde and Brown (1970), it follows that for  any $\delta \in ]0,1]$, 
\begin{multline}
\sup_{x \in {\mathbb R}} \big | {\mathbb P} \big (  \sum_{i=1}^{m}d_{n,m,i}^{\ast  } \leq x \big ) - {\mathbb P} (  s_n N \leq x )\big | \\
\leq K_{\delta} \Big \{s_{n,m} ^{-2-2\delta} \Big ( \mathbb{E}\big(|d_{n,m,i}^{\ast }|^{2+2\delta }) + \mathbb{E}\Big |\sum_{i=1}^{m}\mathbb{E}_{n,m,i-1}\big((d_{n,m,i}^{\ast })^{2}\big)-s_{n,m}^2\Big |%
^{1+\delta }\Big ) \Big \}^{1/(3+2\delta)} \, , \label{ThmHB}
\end{multline}
where $K_{\delta}$ is a positive constant depending only on $\delta$. Assume now that we can prove that there exists  a $\delta$ in $]0,1]$ such that for any positive reals $\varepsilon$ and $M$, 
\begin{equation}
\limsup_{m\rightarrow \infty }\limsup_{n\rightarrow \infty }\sum_{i=1}^{m}%
\mathbb{E}\big(|d_{n,m,i}^{\ast }|^{2+2\delta }\big) \leq u(\varepsilon) \, ,  \label{p4thmgap1}
\end{equation}%
and
\begin{equation}
\limsup_{m\rightarrow \infty }\limsup_{n\rightarrow \infty }\mathbb{E}\Big
|\sum_{i=1}^{m}\mathbb{E}_{n,m,i-1}\big((d_{n,m,i}^{\ast })^{2}\big)-1\Big |%
^{1+\delta }\leq v(M) \, ,  \label{p5thmgap1}
\end{equation}
where $u(\cdot)$ and $v(\cdot)$ are positive functions defined on ${\mathbb R}^+$ such that $v(\cdot)$ does not depend on $\varepsilon$, $\lim_{x \rightarrow 0} u(x)=0$ and $\lim_{x \rightarrow \infty} v(x)=0$. Then starting from \eqref{ThmHB} and noticing that \eqref{p5thmgap1} also implies that for any positive reals $\varepsilon$ and $M$, 
and 
\begin{equation}
\limsup_{m\rightarrow \infty }\limsup_{n\rightarrow \infty }|s_{n,m}^2 -1 |^{1+\delta} \leq v(M) \, ,
\end{equation}
we infer that \eqref{p3thmgap1} will hold. Indeed, by standard arguments, we will get \eqref{p3thmgap1} for every continuous function $f$ with compact support and then \eqref{p3thmgap1} for every continuous and bounded function $f$ by using the fact that every probability measure is tight. Hence, to end the proof of the theorem, it remains to prove that \eqref{p4thmgap1} and \eqref{p5thmgap1} hold. With this aim, we first notice that
$$
\sum_{i=1}^{m}\mathbb{E}\big(|d_{i}^{\ast }|^{2+2\delta }\big) \leq 4(2\varepsilon
)^{2\delta }\sum_{i=1}^{m}\mathbb{E}(d_{i}^{2})
 \leq 4(2\varepsilon )^{2\delta }\Big (1+\mathbb{E}\Big |\sum_{i=1}^{m}
\mathrm{Var}(U_{i}|\mathcal{G}_{i-1})-1\Big |\Big )\,.
$$
Hence, using condition (\ref{C2ap}),  (\ref{p4thmgap1}) follows with $u(\varepsilon) = 4(2\varepsilon )^{2\delta }$. It remains to
prove (\ref{p5thmgap1}). With this aim, using the convexity inequality: $(a+b)^p \leq 2^{p-1} ( a^b +b^p)$ ($p\geq 1$, $a>0$ and $b>0$), 
we first write that   
\begin{equation}
\mathbb{E}\Big |\sum_{i=1}^{m} \mathbb{E}_{i-1}((d_{i}^{\ast })^{2})-1\Big |^{1+\delta }   \label{p6thmgap1}
 \leq 2^{\delta }\mathbb{E}\Big |\sum_{i=1}^{m}\mathbb{E}
_{i-1}(Y_{i}^{2})-1\Big |^{1+\delta }+2^{\delta }\mathbb{E}\Big |
\sum_{i=1}^{m}\big (\mathbb{E}_{i-1}(Y_{i})\big )^{2}\Big |^{1+\delta }\,.
\end{equation}
Now, since $V_{n,m,i}$ is $\mathcal{G}_{n,m,i-1}$-measurable and $\mathbb{E}%
_{n,m,i-1}(d_{n,m,i})=0$ a.s., we infer that
\begin{align*}
\mathbb{E}\Big ( \sum_{i=1}^{m}\big (\mathbb{E}_{i-1}(Y_{i})\big )^{2}\Big )%
^{1+\delta }   
 &\leq  \mathbb{E}\Big (\Big (\sum_{i=1}^{m}\big (\mathbb{E}_{i-1}(d_{i}%
\mathbf{1}_{|d_{i}|>\varepsilon })\big )^{2}\Big )\Big (\sum_{k=1}^{m}%
\mathbf{1}_{V_{k}\leq M}\mathbb{E}_{k-1}(d_{k}^{2})\Big )^{\delta }\Big )\\
 &\leq   M^{\delta }\sum_{i=1}^{m}\mathbb{E}(d_{i}^{2}\mathbf{1}%
_{|d_{i}|>\varepsilon })\,.
\end{align*}
Hence by Lemma \ref{curiouslma},  
\begin{equation}  \label{p7thmgap1} 
\mathbb{E}\Big ( \sum_{i=1}^{m}\big (\mathbb{E}_{i-1}(Y_{i})\big )^{2}\Big )%
^{1+\delta } \leq 12  M^{\delta }\sum_{i=1}^{m}\mathbb{E}(U_{i}^{2}\mathbf{1}%
_{|U_{i}|>\varepsilon/4 })\,.
\end{equation}
On the other hand, using again the fact that $V_{n,m,i}$ is $\mathcal{G}%
_{n,m,i-1}$-measurable and also that $V_{n,m,i}\leq V_{n,m,i+1}$, we derive
that
\begin{eqnarray*}
 \mathbb{E}\Big |\sum_{i=1}^{m}\mathbb{E}_{i-1}(Y_{i}^{2})-1\Big |%
^{1+\delta }
& \leq & \mathbb{E}\Big (\Big (1+\sum_{i=1}^{m}\mathbf{1}_{V_{i}\leq M}\mathbb{%
E}_{i-1}(d_{i}^{2})\Big )^{\delta }\Big |1-\sum_{k=1}^{m}\mathbf{1}%
_{V_{k}\leq M}\mathbb{E}_{k-1}(d_{k}^{2}\mathbf{1}_{|d_{k}|\leq \varepsilon
})\Big |\Big ) \\
& \leq & (M+1)^{\delta }\sum_{i=1}^{m}\mathbb{E}(d_{i}^{2}\mathbf{1}%
_{|d_{i}|>\varepsilon })+(M+1)^{\delta }\mathbb{E}\Big |1-\sum_{k=1}^{m}%
\mathbf{1}_{V_{k}\leq M}\mathbb{E}_{k-1}(d_{k}^{2})\Big | \\
& \leq & (M+1)^{\delta }\sum_{i=1}^{m}\mathbb{E}(d_{i}^{2}\mathbf{1}%
_{|d_{i}|>\varepsilon })+(M+1)^{\delta }\mathbb{E}\Big |\sum_{k=1}^{m}%
\mathbb{E}_{k-1}(d_{k}^{2})-1\Big | \\
& & +(M+1)^{\delta }\mathbb{E}\Big |\mathbf{1}%
_{V_{m}>M}\sum_{k=1}^{m}\mathbb{E}_{k-1}(d_{k}^{2})\Big |\,.
\end{eqnarray*}%
Therefore,
\begin{multline*}
 \mathbb{E}\Big |\sum_{i=1}^{m}\mathbb{E}_{i-1}(Y_{i}^{2})-1\Big |%
^{1+\delta }  
 \leq (M+1)^{\delta }\sum_{i=1}^{m}\mathbb{E}(d_{i}^{2}\mathbf{1}%
_{|d_{i}|>\varepsilon })+2(M+1)^{\delta }\mathbb{E}\Big |\sum_{k=1}^{m}%
\mathbb{E}_{k-1}(d_{k}^{2})-1\Big |   \\
  +\frac{(M+1)^{\delta }}{M}\Big(1+\mathbb{E}\Big |%
\sum_{k=1}^{m}\mathbb{E}_{k-1}(d_{k}^{2})-1\Big |\Big )\, ,
\end{multline*}
which together with Lemma \ref{curiouslma}  
and the fact that $\mathbb{E}_{k-1}(d_{k}^{2})=\mathrm{Var}(U_{k}|\mathcal{G}_{k-1})$ imply that 
\begin{multline}
 \mathbb{E}\Big |\sum_{i=1}^{m}\mathbb{E}_{i-1}(Y_{i}^{2})-1\Big |%
^{1+\delta }   \label{p8thmgap1}
 \leq 12 (M+1)^{\delta }\sum_{i=1}^{m}\mathbb{E}(U_{i}^{2}\mathbf{1}%
_{|U_{i}|>\varepsilon/4 })+2(M+1)^{\delta }\mathbb{E}\Big |\sum_{k=1}^{m}%
\mathrm{Var}(U_{k}|\mathcal{G}_{k-1})-1\Big |   \\
  +\frac{(M+1)^{\delta }}{M}\Big(1+\mathbb{E}\Big |%
\sum_{k=1}^{m}\mathrm{Var}(U_{k}|\mathcal{G}_{k-1})-1\Big |\Big )\, .
\end{multline}
Starting from (\ref{p6thmgap1}) and considering the bounds (\ref{p7thmgap1})
and (\ref{p8thmgap1}) together with
the conditions  (\ref{C2ap}) and (\ref{Lap}), we then infer that (\ref%
{p5thmgap1}) holds for any $\delta \in ]0,1[$ with $v(M)=M^{-1} (M+1)^{\delta}$. This ends the proof of (\ref%
{p3thmgap1}) and then of the theorem. $\square $

\begin{lma}
\label{LGNmart} Assume that
$(d_{n,m,i})_{i\geq 1}$ is an array of random variables in ${\mathbb{L}}^{2}$ adapted to an array $(\mathcal{G}_{n,m,i})_{i\geq 1}$ of nested
sigma fields, and such that for any $i\geq 1$, $\mathbb{E}_{n,m,{i-1}%
}(d_{n,m,i})=0$ almost surely. Suppose that
\begin{equation}
\lim_{m\rightarrow \infty }\limsup_{n\rightarrow \infty }\sum_{i=1}^{m}%
\mathbb{E}(|d_{n,m,i}|\mathbf{1}_{|d_{n,m,i}|>\varepsilon })=0\quad \text{and} \quad
\sum_{i=1}^{m}\mathbb{E}|d_{n,m,i}|<K  \label{1}
\end{equation}%
for some positive constant $K$.
Then
\begin{equation*}
\lim_{m\rightarrow \infty }\limsup_{n\rightarrow \infty }\mathbb{E}\Big |%
\sum_{i=1}^{m}d_{n,m,i}\Big |=0\,.
\end{equation*}
\end{lma}

\noindent \textbf{Proof of Lemma \ref{LGNmart}}. Let $\varepsilon >0$, and
let for any $i\geq 1$,
\begin{equation*}
d_{n,m,i}^{\prime }=d_{n,m,i}\mathbf{1}_{|d_{n,m,i}|\leq \varepsilon }\quad \text{
and }\quad d_{n,m,i}^{\prime \prime }=d_{n,m,i}\mathbf{1}_{|d_{n,m,i}|>\varepsilon
}\,.
\end{equation*}%
With this notation and since $\mathbb{E}_{n,m,{i-1}}(d_{n,m,i})=0$ almost surely,
\begin{equation*}
\sum_{i=1}^{m}d_{n,m,i}=\sum_{i=1}^{m}\big (d_{n,m,i}^{\prime }-\mathbb{E}%
_{n,m,i-1}(d_{n,m,i}^{\prime })\big )+\sum_{i=1}^{m}\big (d_{n,m,i}^{\prime
\prime }-\mathbb{E}_{n,m,i-1}(d_{n,m,i}^{\prime \prime })\big )\,.
\end{equation*}%
Since
\begin{equation*}
\mathbb{E}\big |\big (d_{n,m,i}^{\prime \prime }-\mathbb{E}%
_{n,m,i-1}(d_{n,m,i}^{\prime \prime })\big )\big |\leq 2\mathbb{E}%
(|d_{n,m,i}|\mathbf{1}_{|d_{n,m,i}|>\varepsilon })\,,
\end{equation*}%
by using the first part of (\ref{1}), the lemma will follow if we can prove
that
\begin{equation}
\lim_{\varepsilon \rightarrow 0}\limsup_{m\rightarrow \infty
}\limsup_{n\rightarrow \infty }\mathbb{E}\Big |\sum_{i=1}^{m}\big (%
d_{n,m,i}^{\prime }-\mathbb{E}_{n,m,i-1}(d_{n,m,i}^{\prime })\big )\Big |%
=0\,.  \label{p1LGN}
\end{equation}%
With this aim, it suffices to notice that
$$
\mathbb{E}\Big (\sum_{i=1}^{m} \big (d_{n,m,i}^{\prime }-\mathbb{E}%
_{n,m,i-1}(d_{n,m,i}^{\prime })\big )\Big )^{2}\leq \sum_{i=1}^{m}\mathbb{E}%
(d_{n,m,i}^{\prime })^{2}
 \leq \varepsilon \sum_{i=1}^{m}\mathbb{E}|d_{n,m,i}|\,,
$$
showing that (\ref{p1LGN}) holds under (\ref{1}). $\square $

\begin{lma}
\label{curiouslma} Let $X$ be a real random variable and $\mathcal{F}$ a
sigma-field. For any $p\geq 1$ and any $\varepsilon >0$,
\begin{equation}
\mathbb{E}\big (|X|^{p}\mathbf{1}_{|\mathbb{E}(X|\mathcal{F})|>2\varepsilon }%
\big )\leq 2\,\mathbb{E}\big (|X|^{p}\mathbf{1}_{|X|>\varepsilon }\big )\,,
\label{toapplylind2}
\end{equation}%
and setting $Y=X-\mathbb{E}(X|\mathcal{F})$,
\begin{equation}
\mathbb{E}\big (|X|^{p}\mathbf{1}_{|Y|>3\varepsilon }\big )\leq 2\,\mathbb{E}%
\big (|X|^{p}\mathbf{1}_{|X|>\varepsilon }\big )\quad \text{ and }\quad \mathbb{E}%
(|Y|^{p}\mathbf{1}_{|Y|>4\varepsilon })\leq 3\times 2^{p}\mathbb{E}\big (%
|X|^{p}\mathbf{1}_{|X|>\varepsilon })\,.  \label{toapplylind2bis}
\end{equation}
\end{lma}

\noindent \textbf{Proof of Lemma \ref{curiouslma}.} We first write that
\begin{equation}
|X|^{p}\mathbf{1}_{|\mathbb{E}(X|\mathcal{F})|>2\varepsilon }\leq |X|^{p}%
\mathbf{1}_{|X|>\varepsilon }+\varepsilon ^{p}\mathbf{1}_{|\mathbb{E}(X|%
\mathcal{F})|>2\varepsilon }\,.  \label{p0curiouslma}
\end{equation}%
Notice now that $\{|\mathbb{E}(X|\mathcal{F})|>2\varepsilon \}\subseteq
\{|\mathbb{E}(X\mathbf{1}_{|X|>\varepsilon} )|\mathcal{F})|>\varepsilon \}$,
implying that
\begin{equation}
\varepsilon ^{p}\mathbf{1}_{|\mathbb{E}(X|\mathcal{F})|>2\varepsilon }\leq |%
\mathbb{E}(X\mathbf{1}_{|X|>\varepsilon }|\mathcal{F})|^{p}\leq \mathbb{E}%
\big (|X|^{p}\mathbf{1}_{|X|>\varepsilon }|\mathcal{F}\big )\,,
\label{p2curiouslma}
\end{equation}%
Starting from (\ref{p0curiouslma}), using (\ref{p2curiouslma}) and taking
the expectation, (\ref{toapplylind2}) follows.  
To prove the first part of (\ref%
{toapplylind2bis}),  we start by writing that
\begin{equation*}
|X|^{p}\mathbf{1}_{|Y|>3\varepsilon }\leq |X|^{p}\mathbf{1}_{|X|>\varepsilon
}+\varepsilon ^{p}\mathbf{1}_{|\mathbb{E}(X|\mathcal{F})|>2\varepsilon }\,,
\end{equation*}%
and we use (\ref{p2curiouslma}). To prove the second part of (\ref%
{toapplylind2bis}), we first notice that for any positive reals $a$, $b$ and $\varepsilon $, $(a+b)^{p}%
\mathbf{1}_{a+b>4\varepsilon }\leq 2^{p}a^{p}\mathbf{1}_{a>2\varepsilon
}+2^{p}b^{p}\mathbf{1}_{b>2\varepsilon }$. Therefore
\begin{equation*}
\mathbb{E}(|Y|^{p}\mathbf{1}_{|Y|>4\varepsilon })\leq 2^{p}\mathbb{E}%
(|X|^{p}\mathbf{1}_{|X|>2\varepsilon })+2^{p}\mathbb{E}(|X|^{p}\mathbf{1}%
_{|\mathbb{E}(X|\mathcal{F})|>2\varepsilon })\,. 
\end{equation*}
The second part of (\ref{toapplylind2bis}) then follows by using (\ref{toapplylind2}). $\square$

\section{Ergodic theorem}

We gather below the ergodic theorems used in this paper.
We keep the notations of Section \ref{results} and we define ${\mathcal F}_{ \infty} = \bigvee_{k \in {\mathbb Z}} {\mathcal F}_{k}$.

\begin{lma}
\label{ergodic} Let $Z$ be a ${\mathcal F}_{\infty}$-measurable real-valued random variable in ${\mathbb{L}}^1$. Define $Z_k = Z \circ T^k$ for any $k$ in ${\mathbb{Z}}$. Then
\begin{equation}
\frac{1}{n}\sum_{i=1}^{n}\mathbb{E}_{0}(Z_{i})\rightarrow \mathbb{E}(Z|
\mathcal{I)} \quad  \text{almost surely and in ${\mathbb{L}}^1$,}  \label{ERG}
\end{equation}
and
\begin{equation}
\frac{1}{n}\mathbb{E}_{0}\Big( \max_{1 \leq i \leq n}|Z_{i}|\Big)\rightarrow 0 \quad
\text{almost surely and in ${\mathbb{L}}^1$.}  \label{maxERG}
\end{equation}
\end{lma}

\noindent\textbf{Proof}. By definition of the operator $K$ (see the beginning of Section \ref{results}),
$$
\frac{1}{n}\sum_{i=1}^{n}\mathbb{E}_{0}(Z_{i})= \frac{1}{n} \sum_{i=1}^n K^i(Z)\, .
$$
 Applying the Dunford-Schwartz ergodic theorem (see for
instance Krengel (1985)) we obtain that  $(K(Z)+\cdots +K^n(Z))/n$ converges almost surely
 and in ${\mathbb{L}}^1$ to some $g\in {\mathbb{L}}^1$. We prove now that $g = \mathbb{E}(Z|
\mathcal{I})$. Let $N \in {\mathbb N}$. Define $Z_{0,N} = \mathbb{E}( Z | {\mathcal F}_N)$ and $Z_{k,N}=Z_{0,N} \circ T^k$ for any $k$ in ${\mathbb{Z}}$. From the stationarity of the sequence $(Z_{k,N})_{k \in {\mathbb Z}}$ and the invariance of $\mathbb{E}( Z_{0,N} | \mathcal{I} )$, we have
$$
\Big \Vert  \mathbb{E}( Z_{0,N} | \mathcal{I} ) - \frac{1}{n} \sum_{k=1}^n Z_{k,N} \Big \Vert_1 = \Big \Vert  \mathbb{E}( Z_{0,N} | \mathcal{I} ) - \frac{1}{n} \sum_{k=1 - (n+N)}^{-N} Z_{k,N} \Big \Vert_1 \, .
$$
Both this equality and the ${\mathbb L}^1$-ergodic theorem imply that $\mathbb{E}( Z_{0,N} | \mathcal{I} )$ is the limit in ${\mathbb L}^1$ of a sequence of ${\mathcal F}_0$-measurable random variables. Hence  $\E ( \mathbb{E}( Z_{0,N} | \mathcal{I} ) | {\mathcal F}_0)= \mathbb{E}( Z_{0,N} | \mathcal{I} )$ almost surely. Therefore, noticing that for any $i \in {\mathbb N}$, $\mathbb{E}_{0}(Z_{i}) =\mathbb{E}_{0}(Z_{i,N}) $ and using, once again, the ${\mathbb L}^1$-ergodic theorem, we derive that 
$$
\limsup_{n \rightarrow \infty}\Big \Vert  \frac{1}{n} \sum_{i=1}^{n}\mathbb{E}_{0}(Z_{i}) -\mathbb{E}( Z_{0,N} | \mathcal{I} ) \Big \Vert_1 
 \leq \limsup_{n \rightarrow \infty}\Big \Vert  \frac{1}{n}\sum_{i=1}^{n} Z_{i,N} -\mathbb{E}( Z_{0,N} | \mathcal{I} ) \Big \Vert_1 =0 \, .
$$
Hence the proof will be complete if we show that $ \lim_{N \rightarrow \infty } \Vert \mathbb{E}( Z_{0,N} | \mathcal{I} ) - \mathbb{E}( Z| \mathcal{I} )  \Vert_1=0$. Notice that
$$
\Vert \mathbb{E}( Z_{0,N} | \mathcal{I} ) - \mathbb{E}( Z | \mathcal{I} )  \Vert_1 \leq \Vert \mathbb{E}( Z  | \mathcal{F}_N ) - Z  \Vert_1  \, .
$$
Therefore since $(\mathbb{E}( Z | \mathcal{F}_N ) )_{N \geq 1}$ is an uniformly integrable martingale, and $Z$ is ${\mathcal F}_{\infty}-$measurable, the desired convergence follows by the martingale convergence theorem.

We turn now to the proof of (\ref{maxERG}). With this aim, we notice that
for any $N>0$,
\begin{equation*}
\frac 1 n  \mathbb{E}_0 \Big( \max_{1 \leq i \leq n} |Z_i| \Big) \leq \frac{N}{n} +
\frac{1}{n} \sum_{i=1}^n \mathbb{E}_0 ( |Z_i| \mathbf{1}_{|Z_i| >N}) \, .
\end{equation*}
By using (\ref{ERG}), $ n^{-1}\sum_{i=1}^n
\mathbb{E}_0 ( |Z_i| \mathbf{1}_{|Z_i| >N})$ converges to $ \mathbb{E} ( |Z| \mathbf{1}%
_{|Z| >N} | {\mathcal{I}})$ almost surely and in ${\mathbb{L}}^1$, as $n$ tends to infinity. Therefore
\begin{equation*}
\lim_{N \rightarrow \infty} \limsup_{n \rightarrow \infty}
\frac 1 n \sum_{i=1}^n \mathbb{E}_0 ( |Z_i| \mathbf{1}_{|Z_i| >N}) = 0 \quad  \text{%
almost surely  and in ${\mathbb{L}}^1$,}
\end{equation*}
which ends the proof of (\ref{maxERG}). $\square$

\section{Appendix} 
This section is devoted to the proof of Proposition \ref{propDV}. We shall see that it follows from a slight modification of the example given in Section 5.4 in  Durieu and Voln\'y (2008).

\medskip

We consider the ergodic dynamical system $(\Omega , {\mathcal A}, \mu , T)$ with positive entropy, the sequence $(e_i)_{i \in \mathbb Z}$ of independent identically distributed (i.i.d.) Rademacher random variables with parameter $1/2$, and the $\sigma$-algebra ${\mathcal F}_0$ as described at the beginning of Section 4.1 in Durieu and Voln\'y (2008).  Now,  for any positive integer $k$, we define 
\begin{equation} \label{defNrho}
N_{k} = 4^{k} \ , \ \rho_{k} = \frac{1}{4^{k}} \ , \ \theta_{k} = \frac{1}{{k} 2^{k}} \ , \ \varepsilon_{k} = \frac{1}{{k}^2 4^{3{k}}}\, ,
\end{equation}
and we consider mutually disjoint sets $(A_k)_{k \in {\mathbb Z}}$  by using their Lemma 2 with $2N_k$ instead of $N_k$, and the sequences $(\rho_k)$ and $(\varepsilon_k)$ defined above. In addition to being disjoint, the sets $(A_k)_{k \in {\mathbb N}^*}$ are such that 
\begin{enumerate}
\item[(i)] $\frac{2}{3}\rho_k \leq \mu (A_k ) \leq \rho_k \text{ for all }k \in {\mathbb N}^*$;
\item[(ii)]  for all $k \in {\mathbb N}^*$ and all $i,j \in \{ 0, \dots, 2 N_k \}$, $\mu ( T^{-i} A_k \Delta T^{-j} A_k ) \leq \varepsilon_k$, 
\end{enumerate}
The function $f$ is then defined as 
\begin{equation} \label{deffunction}
\text{$f = \sum_{k \geq 1} f_k {\bf 1}_{A_k}$    with $f_k = \theta_k \sum_{j=N_k +1}^{2 N_k} e_{-j}$}\, .
\end{equation}
The function $f$ defined in \eqref{deffunction} is centered, ${\mathcal F}_0$-measurable and, since $\sum_{k\geq 1} \theta_k^2 N_k \rho_k < \infty$, it belongs to ${\mathbb L}^2$ (see Proposition 7 in Durieu and Voln\'y (2008)). 

Let now $X_i=f \circ T^i$ for any $i \in {\mathbb Z}$. This sequence is adapted to the stationary and nondecreasing sequence of $\sigma$-algebras $({\mathcal F}_i )_{i \in {\mathbb Z}}$ where ${\mathcal F}_i = T^{-i}({\mathcal F}_0)$. Let us first prove that the sequence $(X_i)_{i \in {\mathbb Z}}$  satisfies the condition \eqref{cond}. With this aim, we first  emphasize some additional important properties of $(e_i)_{i \in \mathbb Z}$ and of $(A_k)_{k \in {\mathbb Z}}$. First, the sequence $(e_i)_{i \in \mathbb Z}$ is adapted to $({\mathcal F}_i)_{i \in {\mathbb Z}}$ and ${\mathbb E} (e_i |{\mathcal F}_0) = e_i {\bf 1}_{i \leq 0}$ almost surely. Second, for all $k$ and $i$, ${\bf 1}_{A_k} \circ T^i$ is ${\mathcal F}_0$-measurable. Finally, the $e_i$'s and the ${\bf 1}_{A_k}$'s are independent for all $i$ and $k$. As in relation (4) in Durieu and Voln\'y (2008), we then write that for any $i \in {\mathbb N}$,
\begin{align} \label{consta1CE}
{\mathbb E} (X_i |{\mathcal F}_0) & = \sum_{k \geq 1} \theta_k\sum_{j = N_k+1}^{2N_k} e_{i-j}{\bf 1}_{T^{-i}(A_k)}   {\bf 1}_{i \leq j} \nonumber \\
& = \sum_{k \geq 1}\theta_k \sum_{j = N_k+1}^{2N_k} e_{i-j}{\bf 1}_{A_k}   {\bf 1}_{i \leq j} +\sum_{k \geq 1} \theta_k \sum_{j = N_k+1}^{2N_k} e_{i-j} ( {\bf 1}_{T^{-i}(A_k) \backslash   A_k}  -{\bf 1}_{A_k \backslash  T^{-i}(A_k)  }  )  {\bf 1}_{i \leq j} \, .
\end{align}
Using Item (ii) above, and the fact that the $e_j$'s are bounded by one, we obtain
\begin{align} \label{consta2CE}
\sum_{i \geq 0} \Big \Vert \sum_{k \geq 1} & \theta_k \sum_{j = N_k+1}^{2N_k } e_{i-j} {\bf 1}_{i \leq j}  ( {\bf 1}_{T^{-i}(A_k) \backslash   A_k}  -{\bf 1}_{A_k \backslash  T^{-i}(A_k)  }  )  \Big \Vert_2  \nonumber \\
& \leq \sum_{k \geq 1} \sum_{i =0}^{N_k} \theta_k N_k \big  ( \mu ( T^{-i}(A_k) \Delta  A_k) \big )^{1/2} +\sum_{k \geq 1} \theta_k \sum_{i =N_k+1}^{2N_k} \sum_{j =i}^{2N_k } \big  ( \mu ( T^{-i}(A_k) \Delta  A_k) \big )^{1/2}  \nonumber  \\
& \leq 2 \sum_{k \geq 1} \theta_k N_k (N_k +1) \sqrt{\varepsilon_k } \, .
\end{align}
Since, by \eqref{defNrho}, $\sum_{k \geq 1} \theta_k N_k^2 \sqrt{\varepsilon_k } < \infty$, in order to prove that \eqref{cond} holds, it is enough to show that 
\begin{equation} \label{CEbut1}
\sum_{i \geq 0} \Big \Vert f \sum_{k \geq 1}\theta_k \sum_{j = N_k+1}^{2N_k} e_{i-j}{\bf 1}_{A_k}   {\bf 1}_{i \leq j}\Big \Vert_1 < \infty  \, .
\end{equation}
By disjointness of the $A_k$'s,
\begin{multline} \label{85bis}
\sum_{i \geq 0}  \Big \Vert f \sum_{k \geq 1}\theta_k \sum_{j = N_k+1}^{2N_k} e_{i-j}{\bf 1}_{A_k}  {\bf 1}_{i \leq j} \Big \Vert_1 =\sum_{i \geq 0} \Big \Vert  \sum_{k \geq 1}\theta^2_k \Big ( \sum_{j = N_k+1}^{2N_k} e_{i-j}  {\bf 1}_{i \leq j} \Big ) \Big ( \sum_{\ell = N_k+1}^{2N_k} e_{-\ell} \Big )  {\bf 1}_{A_k}  \Big \Vert_1 \\
 =\sum_{i \geq 0} \Big \Vert  \sum_{k \geq 1} {\bf 1}_{i \leq N_k} \theta^2_k \Big ( \sum_{j = N_k+1}^{2N_k} e_{i-j}  \Big ) \Big ( \sum_{\ell = N_k+1}^{2N_k} e_{-\ell} \Big )  {\bf 1}_{A_k}  \Big \Vert_1 \\+ \sum_{i \geq 0} \Big \Vert  \sum_{k \geq 1} {\bf 1}_{N_k +1 \leq i \leq 2 N_k} \theta^2_k \Big ( \sum_{j = i}^{2N_k} e_{i-j}  \Big ) \Big ( \sum_{\ell = N_k+1}^{2N_k} e_{-\ell} \Big )  {\bf 1}_{A_k}  \Big \Vert_1 \, .
\end{multline}
Now, by independence between the  $e_i$'s and the ${\bf 1}_{A_k}$'s,
\begin{multline*}
\sum_{i \geq 0} \Big \Vert  \sum_{k \geq 1} {\bf 1}_{i \leq N_k} \theta^2_k \Big ( \sum_{j = N_k+1}^{2N_k} e_{i-j}  \Big ) \Big ( \sum_{\ell = N_k+1}^{2N_k} e_{-\ell} \Big )  {\bf 1}_{A_k}  \Big \Vert_1 \\
\leq \sum_{i \geq 0}  \sum_{k \geq 1} {\bf 1}_{i \leq N_k} \theta^2_k \Big \Vert \Big ( \sum_{j = N_k+1}^{2N_k} e_{i-j}  \Big ) \Big ( \sum_{\ell = N_k+1}^{2N_k} e_{-\ell} \Big )   \Big \Vert_1 \mu (A_k) \, .
\end{multline*}
Since the  $e_i$'s are i.i.d., centered and with variance one, we have
$$
\Big \Vert \Big ( \sum_{j = N_k+1}^{2N_k} e_{i-j}  \Big ) \Big ( \sum_{\ell = N_k+1}^{2N_k} e_{-\ell} \Big )   \Big \Vert_1  \leq \Big \Vert  \sum_{j = N_k+1}^{2N_k} e_{i-j}  \Big \Vert_2 \Big \Vert \sum_{\ell = N_k+1}^{2N_k} e_{-\ell}    \Big \Vert_2 \leq N_k \, .
$$
The second term in the right-hand side of \eqref{85bis} can be handled similarly. So overall, we infer that
$$
\sum_{i \geq 0}  \Big \Vert f \sum_{k \geq 1}\theta_k \sum_{j = N_k+1}^{2N_k} e_{i-j}{\bf 1}_{A_k}  {\bf 1}_{i \leq j} \Big \Vert_1  \leq 2 \sum_{k \geq 1} \theta^2_k N_k^2 \rho_k \, ,
$$
which is finite according to  \eqref{defNrho}. This ends the proof of \eqref{CEbut1} and then of the fact that the sequence $(X_i)_{i \in {\mathbb Z}}$ satisfies the condition \eqref{cond}.

\medskip

Let us prove now that the condition \eqref{condGordinL1} fails for the sequence $(X_i)_{i \in {\mathbb Z}}$ defined above. With this aim, we shall prove that
\begin{equation} \label{CEgordinl1}
\sup_{K \in {\mathbb N}} \Vert \E_0(S_{N_K} )\Vert_1 = \infty \, .
\end{equation}
Starting from \eqref{consta1CE} and using \eqref{consta2CE},  it suffices to prove that 
\begin{equation} \label{CEgordinl12}
\sup_{K \in {\mathbb N}} \Big \Vert \sum_{i=1}^{N_K} \sum_{k \geq 1} \theta_k \sum_{j=N_k +1}^{2N_k} e_{i-j} {\bf 1}_{i \leq j}{\bf 1}_{A_k}  \Big \Vert_1 = \infty \, .
\end{equation}
Let $K \geq 3$. By disjointness of the $A_k$'s and by independence between the  $e_i$'s and the ${\bf 1}_{A_k}$'s, it follows that 
$$
\Big \Vert \sum_{i=1}^{N_K} \sum_{k \geq 1} \theta_k  \sum_{j=N_k +1}^{2N_k} e_{i-j} {\bf 1}_{i \leq j}{\bf 1}_{A_k} \Big  \Vert_1  = \sum_{k \geq 1} \theta_k \mu (A_k) \E \Big |  \sum_{j=N_k +1}^{2N_k} \sum_{i=1}^{N_K} e_{i-j} {\bf 1}_{i \leq j} \Big | \, .
$$
Therefore
$$
\Big \Vert \sum_{i=1}^{N_K} \sum_{k \geq 1} \theta_k  \sum_{j=N_k +1}^{2N_k} e_{i-j} {\bf 1}_{i \leq j}{\bf 1}_{A_k} \Big \Vert_1  \geq  \sum_{k = 1}^{K-1} \theta_k \mu (A_k) \E \Big |  \sum_{\ell=(N_k +1 -N_K)\vee 0}^{2N_k-1} e_{- \ell} \sum_{i=1}^{N_K} {\bf 1}_{i \geq N_k +1 - \ell} {\bf 1}_{i \leq 2N_k - \ell}  \Big | \, .
$$
Notice now that for any $k \in \{1, \dots, K-1\}$, $
N_k +1-N_K \leq N_{K-1} +1 -N_K \leq 0 $ and for any $\ell \geq 0$, $2N_k - \ell \leq N_K$ (since $2N_k -N_K \leq 2N_{K-1} -N_K \leq 0$). Hence  
$$
\Big \Vert \sum_{i=1}^{N_K} \sum_{k \geq 1} \theta_k  \sum_{j=N_k +1}^{2N_k} e_{i-j} {\bf 1}_{i \leq j}{\bf 1}_{A_k} \Big  \Vert_1  \geq  \sum_{k = 1}^{K-1} \theta_k \mu (A_k) \E \Big | N_k \sum_{\ell= 0}^{N_k} e_{-\ell}  +\sum_{\ell=N_k +1}^{2N_k-1} ( 2N_k - \ell ) e_{-\ell}  \Big | \, .
$$
Next, by using the Marcinkiewicz-Zygmund's inequality together with Item (i) above, we get that there exists a positive constant $A$ such that
\begin{align*}
\Big \Vert \sum_{i=1}^{N_K} \sum_{k \geq 1} \theta_k  & \sum_{j=N_k +1}^{2N_k} e_{i-j} {\bf 1}_{i \leq j}{\bf 1}_{A_k} \Big \Vert_1  \geq  A \sum_{k = 1}^{K-1} \theta_k \rho_k \E \Big ( N^2_k \sum_{\ell= 0}^{N_k} e^2_{-\ell}  +\sum_{\ell=N_k +1}^{2N_k-1} ( 2N_k - \ell )^2 e^2_{-\ell}  \Big )^{1/2} \\
& \geq   A \sum_{k = 1}^{K-1} \theta_k \rho_k N^{3/2}_k \geq A \ln (K-1) \, ,
\end{align*} 
which proves \eqref{CEgordinl12} and therefore \eqref{CEgordinl1}.

\medskip

Let us prove now that the condition \eqref{condMW*} fails for the sequence $(X_i)_{i \in {\mathbb Z}}$ defined above. Following the computations page 339 in Durieu and Voln\'y (2006),  it suffices to prove that 
\begin{equation} \label{CEMW1}
\sum_{n \geq 1} \frac{1}{n^{3/2}} \Big ( \sum_{k \geq 1}{\bf 1}_{2N_k \leq n } \theta_k^2 N_k^3 \rho_k  \Big )^{1/2} =\infty \, .
\end{equation}
Since $N_k = 4^k$, using \eqref{defNrho}, we get that 
$$
\sum_{k \geq 1}{\bf 1}_{2N_k \leq n } \theta_k^2 N_k^3 \rho_k = \sum_{k = 1}^{[(\ln n - \ln 2) /(2 \ln 2) ]} \frac{4^k}{k^2} \geq C \frac{n}{( \ln n)^2} \, ,
$$
where $C$ is a positive constant. This shows \eqref{CEMW1} and then that \eqref{condMW*} fails.

\medskip

Let us prove now that the Hannan-Heyde condition \eqref{HHcond} fails for the sequence $(X_i)_{i \in {\mathbb Z}}$ defined above. With this aim, we first notice that 
$$
 {\mathbb E} (X_i |{\mathcal F}_0) - {\mathbb E} (X_i |{\mathcal F}_{-1}) =e_{0} \sum_{k \geq 1} \theta_k\sum_{j = N_k+1}^{2N_k} {\bf 1}_{T^{-i}(A_k)}   {\bf 1}_{i =j} \, .$$
Proceeding as in \eqref{consta1CE} and since $\sum_{k \geq 1} \theta_k N_k \sqrt{\varepsilon_k } < \infty$, it suffices to prove that 
\begin{equation} \label{CEbut2}
\sum_{i \geq 1} \Big \Vert  \sum_{k \geq 1}\theta_k \sum_{j = N_k+1}^{2N_k} {\bf 1}_{i = j} {\bf 1}_{A_k}   \Big \Vert_2 =\infty  \, .
\end{equation}
But, by disjointness of the $A_k$'s,
$$
\sum_{i \geq 1} \Big \Vert  \sum_{k \geq 1}\theta_k \sum_{j = N_k+1}^{2N_k} {\bf 1}_{i = j} {\bf 1}_{A_k}   \Big \Vert_2 
 =\sum_{i \geq 1} \Big ( \sum_{k \geq 1} \theta^2_k \Big ( \sum_{j = N_k+1}^{2N_k}  {\bf 1}_{i = j}  \Big )^2 \mu (A_k) \Big )^{1/2} \, .
$$
Therefore, since $N_k=2^{2k}$,
\begin{align*}
& \sum_{i \geq 1}   \Big \Vert  \sum_{k \geq 1}\theta_k \sum_{j = N_k+1}^{2N_k} {\bf 1}_{i = j} {\bf 1}_{A_k}   \Big \Vert_2 
= \sum_{\ell \geq 0} \sum_{i =2^{2 \ell} +1 }^{2^{2(\ell +1)}} \Big ( \sum_{k \geq 1} \theta^2_k \Big ( \sum_{j = 2^{2k} +1}^{2^{2k+1}}  {\bf 1}_{i = j}  \Big )^2 \mu (A_k)  \Big )^{1/2} \\
& \geq  \sum_{\ell \geq 0} \sum_{i =2^{2\ell} +1 }^{2^{2 \ell +1}} \Big ( \sum_{k \geq 1} \theta^2_k \Big ( \sum_{j = 2^{2k} +1}^{2^{2k+1}}  {\bf 1}_{i = j}  \Big )^2 \mu (A_k)  \Big )^{1/2}  = \sum_{\ell \geq 0} \sum_{i =2^{\ell} +1 }^{2^{2\ell +1}} \big (  \theta^2_{\ell} \mu (A_{\ell}) \big )^{1/2} \geq \frac{\sqrt 2}{\sqrt 3}\sum_{\ell \geq 0} 2^{2\ell }  \theta_{\ell} \sqrt{\rho_{\ell}} 
\, ,
\end{align*}
which does not converge according to  \eqref{defNrho}. This ends the proof of \eqref{CEbut2} and then of the fact that the sequence $(X_i)_{i \in {\mathbb Z}}$ does not satisfy the Hannan-Heyde condition \eqref{HHcond}. $\square$

\medskip

\noindent {\bf Acknowledgements.} The authors would like to thank the referee for
carefully reading the manuscript and for numerous suggestions that improved
the presentation of this paper. The authors are also indebted to Christophe Cuny for helpful discussions.

\end{document}